
\magnification=1200

\pageno=1
\baselineskip=18pt

\parindent=0pt

\font\bigbf=cmbx10 scaled\magstep1
\raggedbottom

\input eplain

\newfam\bbbfam
\font\bbbten=msbm10
\font\bbbseven=msbm7
\font\bbbfive=msbm5
\textfont\bbbfam=\bbbten
\scriptfont\bbbfam=\bbbseven
\scriptscriptfont\bbbfam=\bbbfive
\def\bbb{\fam=\bbbfam}

\font\teneufm=eufm10
\font\seveneufm=eufm7
\font\fiveeufm=eufm5
\newfam\eufmfam
\textfont\eufmfam=\teneufm
\scriptfont\eufmfam=\seveneufm
\scriptscriptfont\eufmfam=\fiveeufm
\def\eufm#1{{\fam\eufmfam\relax#1}}

\font\teneufb=eufb10
\font\seveneufb=eufb7
\font\fiveeufb=eufb5
\newfam\eufbfam
\textfont\eufbfam=\teneufb
\scriptfont\eufbfam=\seveneufb
\scriptscriptfont\eufbfam=\fiveeufb

\font\teneurm=eurm10
\font\seveneurm=eurm7
\font\fiveeurm=eurm5
\newfam\eurmfam
\textfont\eurmfam=\teneurm
\scriptfont\eurmfam=\seveneurm
\scriptscriptfont\eurmfam=\fiveeurm

\font\teneurb=eurb10
\font\seveneurb=eurb7
\font\fiveeurb=eurb5
\newfam\eurbfam
\textfont\eurbfam=\teneurb
\scriptfont\eurbfam=\seveneurb
\scriptscriptfont\eurbfam=\fiveeurb

\font\tenscr=rsfs10
\font\sevenscr=rsfs7
\font\fivescr=rsfs5
\skewchar\tenscr='177 \skewchar\sevenscr='177 \skewchar\fivescr='177
\newfam\scrfam \textfont\scrfam=\tenscr \scriptfont\scrfam=\sevenscr
\scriptscriptfont\scrfam=\fivescr
\def\scr#1{{\fam\scrfam#1}}

\rm
\centerline{\bigbf\ Cusp forms of weight 1/2 and pairs of quadratic forms}
 \vskip10pt
\centerline{by Andr\'as BIR\'O\footnote{}{Research partially supported by the NKFIH (National Research, Development and Innovation Office) Grants No. K135885, K143876, and by the R\'enyi Int\'ezet Lend\"ulet Automorphic Research Group }}
 \footnote{}{} \footnote{}{2020

}
\hfill\break
\centerline{HUN-REN Alfr\'ed R\'enyi Institute of Mathematics}

\centerline{ 1053 Budapest, Re\'altanoda u. 13-15., Hungary; e-mail: biro.andras@renyi.hu}
 \vskip20pt

\noindent {\bf Abstract}. We prove a spectral summation
formula for the product of four Fourier coefficients of
half-integral weight cusp forms in Kohnen's subspace.
The other side of the formula involves certain
generalized class numbers of pairs of quadratic forms with
integer coefficients.

\vskip20pt

\noindent{\bf 1. Introduction }
\medskip

{\bf 1.1. Informal discussion of the result.} The famous
Bruggeman-Kuznetsov formula (see e.g. [I], Theorem 9.3)
relates a spectral sum of the product of two Fourier
coefficients of cusp forms of weight $0$ to a sum of
Kloosterman sums. This result has been generalized to arbitrary weights in [P], and to Kohnen's subspace of cusp forms of half-integral weight in [B1], [A-A], [A-D], [B-C].

In the present paper we prove a formula which relates a spectral sum of the product of {\it four} Fourier coefficients of cusp forms in Kohnen's subspace to certain arithmetic objects. The arithmetic objects here are weighted summations over the $SL_2({\bf Z})$-equivalence classes of pairs of integral quadratic forms with given discriminants and codiscriminant. We can call these summations generalized class numbers.

The spectral sum we consider here contains two positive and two negative Fourier coefficients (see Theorem 1.1), but similar formulas can be proved under different assumptions on the signs of the Fourier coefficients by modifying the present proof accordingly. We will discuss in Subsection 1.8 how the present proof can be modified in the other cases. In particular, we prove an integral repesentation for the product of two negative Fourier coefficients of a half-integral weight cusp form in Kohnen's subspace, see Theorem 1.2. This theorem can be used in extending Theorem 1.1 for the more difficult cases when there are more negative than positive Fourier coefficients.

Class numbers of pairs of quadratic forms occurred in the
context of automorphic functions in our recent paper [B2],
where we expressed the inner product of two automorphic
functions by a summation of such class numbers, see [B2,
Lemma 2.2]. A variant of that lemma plays a very important role in the present paper.

A summation formula of a different shape for the product of four half-integral weight coefficients, in the case when two of the factors are first coefficients, was proved in [B-C], see Theorem 3 and the lines below that theorem on p 1329 of [B-C]. That theorem was one of the main ingredients in the proof of the main result of [B-C].

In the following subsections we give the definitions
needed to state Theorem 1.1.

{\bf 1.2. Weight} $ ${\bf 1/2 cusp forms and the Shimura lift.} Let $\bbb
H$ be the open upper half-plane. The elements of the
group $PSL_2({\bf R})$ act on $\bbb H$ by linear fractional transformations, these are isometries of the
hyperbolic plane. Let $d\mu_z={{dxdy}\over {y^2}}$, this measure is invariant with respect to
the action of $PSL_2({\bf R})$ on $\bbb H$.

We write
$$\Gamma_0(4)=\left\{\left(\matrix{a&b\cr
c&d\cr}
\right)\in SL_2({\bf Z}):\hbox{\rm \ $c\equiv 0\hbox{\rm \ (mod 4}
)$}\right\}.$$
Let ${\cal F}_1$ be a fundamental domain of $SL_2({\bf Z})$ in $\bbb
H$, and let
${\cal F}_4$ be a fundamental domain of $\Gamma_0(4)$ in $\bbb H$. Let
us write
$$\left(f_1,f_2\right)_1:=\int_{{\cal F}_1}f_1(z)\overline {f_2(z
)}d\mu_z,\;\left(f_1,f_2\right)_4:=\int_{{\cal F}_4}f_1(z)\overline {
f_2(z)}d\mu_z.$$
For a complex number $z\neq 0$ we set its argument in
$(-\pi ,\pi ]$, and write $\log z=\log\left|z\right|+i\arg z,$ where $\log\left
|z\right|$ is
real. We define the power $z^s$ for any $s\in {\bf C}$ by
$z^s=e^{s\log z}$. We write $e(x)=e^{2\pi ix}$.

For $z\in\bbb H$ we define
$$B_0(z):=\left(\hbox{\rm Im$z$}\right)^{{1\over 4}}\theta\left(z\right
)=\left(\hbox{\rm Im$z$}\right)^{{1\over 4}}\sum_{m=-\infty}^{\infty}
e(m^2z).$$
Then
$$B_0(\gamma z)=\nu (\gamma )\left({{j_{\gamma}(z)}\over {\left|j_{
\gamma}(z)\right|}}\right)^{1/2}B_0(z)\hbox{\rm \ for $\gamma\in\Gamma_
0(4)$}$$
with a well-known multiplier system $\nu$, where for $\gamma =$$\left
(\matrix{a&b\cr
c&d\cr}
\right)\in SL_2({\bf R})$ we write $j_{\gamma}(z)=cz+d$.

The three cusps for $\Gamma_0(4)$ are $\infty$, $0$ and $-{1\over
2}$. If $\eufm{a}$ denotes
one of these cusps, we take a scaling matrix
$\sigma_{\eufm{a}}\in SL_2({\bf R})$ as it is explained on p 42 of [I]. We can easily
see that one can take
$$\sigma_{\infty}:=\left(\matrix{1&0\cr
0&1\cr}
\right),\qquad\sigma_0:=\left(\matrix{0&{{-1}\over 2}\cr
2&0\cr}
\right),\qquad\sigma_{-{1\over 2}}:=\left(\matrix{-1&{{-1}\over 2}\cr
2&0\cr}
\right).\eqno (1.1)$$
We define $\chi_{\eufm{a}}$ by
$$\nu\left(\sigma_{\eufm{a}}\left(\matrix{1&1\cr
0&1\cr}
\right)\sigma_{\eufm{a}}^{-1}\right)=e(-\chi_{\eufm{a}}),\qquad 0
\le\chi_{\eufm{a}}<1.$$
It is easy to check that $\chi_{\infty}=\chi_0=0$, and $\chi_{-{1\over
2}}={3\over 4}$. The
cusps with $\chi_{\eufm{a}}=0$ are said to be singular. So the singular cusps of $
\Gamma_0(4)$ are $\infty$ and $0$.

The only cusp for $SL_2({\bf Z})$ is $\infty$.

Introduce the hyperbolic Laplace operator of weight $k$
for any real $k$:
$$\Delta_k:=y^2\left({{\partial^2}\over {\partial x^2}}+{{\partial^
2}\over {\partial y^2}}\right)-iky{{\partial}\over {\partial x}}.$$
Let $k=0$ or $k={1\over 2}$, and $\Gamma^{(0)}:=SL(2,{\bf Z})$, $
\Gamma^{(1/2)}:=\Gamma_0(4)$. We say that a function $f$ on $\bbb
H$ is
an automorphic function of weight $k$ for $\Gamma^{(k)}$ if it satisfies the transformation formula
$$f(\gamma z)=\left({{j_{\gamma}(z)}\over {\left|j_{\gamma}(z)\right
|}}\right)^kf(z)\cdot\left\{\matrix{1&\quad {\rm i}{\rm f}&k=0\cr
\cr
\nu (\gamma )&\quad {\rm i}{\rm f}&k={1\over 2}\cr}
\right.$$
for any $z\in\bbb H$ and $\gamma\in\Gamma^{(k)}$. The operator $\Delta_
k$ acts on
smooth automorphic functions of weight $k$. We say that a
smooth automorphic function $f$ is a Maass form of
weight $k$ for $\Gamma^{(k)}$ if it has at most polynomial growth at
the cusps of $\Gamma^{(k)}$ and it is an eigenfunction of $\Delta_
k$. If a Maass form $f$ has exponential decay at
all of the cusps of $\Gamma^{(k)}$, it is called a cusp form.

If $f$ is a cusp form of weight $k$ and $\Delta_kf=s(s-1)f$ with some Re$
s\ge{1\over 2}$, $s={1\over 2}+it,$ then one has the Fourier expansion
$$f(z)=\sum_{m\neq 0}\rho_f(m)W_{{k\over 2}{\rm s}{\rm g}{\rm n}(
m),it}\left(4\pi\left|m\right|y\right)e\left(mx\right)\eqno (1.2)$$
for $z=x+iy\in\bbb H$, where $W_{\alpha ,\beta}$ is the Whittaker
function (see [G-R], p 1014). The number $\rho_f(m)$ is called the $
m$th Fourier coefficient
of $f$.

Denote by $L_{1/2}^2(D_4)$ the space of automorphic functions of
weight $1/2$ for $\Gamma_0(4)$ for which we have $(f,f)_4<\infty$. Let $
V$ be the subspace of $L_{1/2}^2$$(D_4)$ consisting of cuspidal
functions $f$, which means that the zeroth Fourier
coefficient of $f$ is $0$ in the two singular cusps of $\Gamma_0(
4)$,
i.e.
$$\int_0^1f\left(\sigma_{\eufm{a}}\left(x+iy\right)\right)\left({{
j_{\sigma_{\eufm{a}}}\left(x+iy\right)}\over {\left|j_{\sigma_{\eufm{
a}}}\left(x+iy\right)\right|}}\right)^{-1/2}dx=0$$
for every $y>0$ and $\eufm{a}=0,\infty$, see (1.1).

Let the operator $L$ have the same meaning as on p. 195 of [K-S]. Let $
V^{+}$ be the subspace of
$V$ with $L$-eigenvalue $1$. This space is called Kohnen's
subspace. It is known that a cusp form $F$ of weight $1/2$ for $\Gamma_
0(4)$
belongs to $V^{+}$ if and only if $\rho_F(m)=0$ for every integer
$m\equiv 2,3(4)$. The holomorphic analogue of this statement is
proved in [K2], Proposition 1, and the proof given there can be
generalized to our case.

The weight $1/2$ Hecke operators $T_{p^2}$:$V^{+}\rightarrow V^{+}$ are
defined in [K-S] for every prime $p>2$. The operators $T_{p^2}$
(for primes $p>2$) form a commouting family of linear self-adjoint operators
$V^{+}\rightarrow V^{+}$, and each of these operators commute with $
\Delta_{1/2}$.

As on p 196 of [K-S], let $F_j$ ($j=1,2,\ldots$) be an orthonormal basis of $
V^{+}$ consisting of
common eigenfunctions of $\Delta_{{1\over 2}}$ and the Hecke operators
$T_{p^2}$ for primes $p>2$. The $F_j$'s are Maass cusp forms of
weight ${1\over 2}$ for the group $\Gamma_0(4)$. Let
$\Delta_{1/2}F_j$$=\left(-{1\over 4}-r_j^2\right)F_j$, where $r_j
\ge 0$ for $j\ge 1$. Denote the Fourier coefficients of $F_j$ by $
b_j(m)$, i.e.
$$b_j(m)=\rho_{F_j}(m).\eqno (1.3)$$
If $j\ge 1$ is an integer, the Shimura lift $\hbox{\rm Shim$F_j$}$ is defined
in [K-S], pp 196-197 under the condition $b_j(1)\neq 0$, and
it is defined also without this condition on p 981 of
[D-I-T]. For every $j\ge 1$ the function $\hbox{\rm Shim$F_j$}$ is a Maass
cusp form of weight $0$ for $SL_2({\rm Z})$, which is a simultaneous Hecke eigenform, even and Hecke
normalized (i. e. for its Fourier coefficients $a(n)$ we
have $a(1)=1$ and $a(n)=a(-n)$). We will discuss the
Shimura lift in more detail in Subsection 3.1.

{\bf 1.3. Zagier} $L${\bf -functions.} If $D$ is a fundamental discriminant, let $
\chi_D$ be the Diriclet
character associated to $D$, it is a real primitive character
of conductor $|D|$. It is given by the symbol $\left({D\over {}}\right
)$, i.e. we
have $\chi_D\left(n\right)=\left({D\over n}\right)$ for every integer $
n$, see [D], p 40.

Let $\zeta$ be the Riemann zeta function. For a Dirichlet character $
\chi$ Let $L\left(s,\chi\right)$ be the Dirichlet
$L$-function associated to $\chi$.

If $\delta$ is a nonzero integer with $\delta\equiv 0,1\hbox{\rm \ (mod 4}
)$, we define
for $\hbox{\rm Re$s>1$}$ the Zagier $L$-series $L\left(s,\delta\right
)$ in the following way:
$$L\left(s,\delta\right)={{\zeta\left(2s\right)}\over {\zeta\left
(s\right)}}\sum_{q=1}^{\infty}{1\over {q^s}}\left(\sum_{r\,{\rm m}
{\rm o}{\rm d}\,2q,\;r^2\equiv\delta\left(4q\right)}1\right).\eqno
(1.4)$$
It is known that if $\delta =Dl^2$ with a fundamental discriminant $
D$ and a positive
integer $l$, then
$$L\left(s,\delta\right)=L\left(s,\chi_D\right)l^{{1\over 2}-s}\sum_{
l_1l_2=l}\chi_D\left(l_1\right){{\mu\left(l_1\right)}\over {\sqrt {
l_1}}}\tau_s\left(l_2\right)\eqno (1.5)$$
with $\tau_s\left(k\right):=k^{s-{1\over 2}}$$\sum_{a|k}a^{1-2s}$, see [S-Y, (4) and (5)].
We see that $L\left(s,\delta\right)$ has a meromorphic continuation to
the complex plane, and if $\delta$ is not a square, then it is an
entire function. We see also that for a fundamental discriminant $
D$ we have
$$\hbox{\rm $L\left(s,D\right)=L\left(s,\chi_D\right)$.}\eqno (1.
6)$$
Let us use the notation
$$L^{\ast}\left(s,\delta\right):=L\left(s,\delta\right)\left|\delta\right
|^{s/2}.\eqno (1.7)$$
{\bf 1.4. Quadratic forms.} If $\delta$ is a nonzero integer with $
\delta\equiv 0,1\hbox{\rm \ (mod$ $ $4$})$, let
$$\scr{Q}_{\delta}:=\left\{Q(X,Y)=AX^2+BXY+CY^2:\;A,B,C\in {\bf Z}
,\;B^2-4AC=\delta\right\}.\eqno (1.8)$$
If $\tau =\left(\matrix{a&b\cr
c&d\cr}
\right)\in SL_2({\bf Z})$ and $Q$ is a quadratic form, let us
define the quadratic form $Q^{\tau}$ by $Q^{\tau}\left(X,Y\right)
=Q\left(aX+bY,cX+dY\right)$. The
group $SL_2({\bf Z})$ acts in this way on $\scr{Q}_{\delta}$.

If $Q(X,Y)=AX^2+BXY+CY^2$ is an element of $Q_{\delta}$ with some
$\delta <0$, let $z_Q$ be the unique root in $\bbb H$ of $Az^2+Bz
+C$, let
$$C(Q)=\left\{\gamma\in PSL_2({\bf Z}):\;\gamma z_Q=z_Q\right\},$$
and $M_Q=\left|C(Q)\right|$.

If $Q(X,Y)=aX^2+bXY+cY^2$ is a quadratic form with
integer coefficients, $d=b^2-4ac$ is its discriminant,
$d\neq 0$ and $D$ is a fundamental discriminant with $D|d$ and
$d/D\equiv 0,1$ (mod $4$), define
$$\omega_D\left(Q\right)=\left\{\matrix{0&\quad {\rm i}{\rm f}&\left
(a,b,c,D\right)>1,\cr
\cr
\left({D\over r}\right)&\quad {\rm i}{\rm f}&\qquad\left(a,b,c,D\right
)=1,\cr}
\right.$$
where $r$ is any number represented by $Q$ with $\left(r,D\right)
=1$.
The symbol $\omega_D\left(Q\right)$ is well-defined, and it depends only
on the $SL_2({\bf Z})$-equivalence class of $Q$ (see [K1]).

For $d_1,d_2,t\in {\bf Z}$, $d_i\equiv 0,1$ (mod $4$) for $i= 1,2$, let ${\cal Q}_{d_1,d_2,t}$ be the subset of $\scr{
Q}_{d_1}\times\scr{Q}_{d_2}$ consisting of those pairs $\left(Q_1
,Q_2\right)$ of quadratic forms having codiscriminant $t$. In other words, writing
$$Q_1\left(X,Y\right)=A_1X^2+B_1XY+C_1Y^2,\;Q_2\left(X,Y\right)=A_
2X^2+B_2XY+C_2Y^2\eqno (1.9)$$
we require that the discriminant of $Q_j$ is $d_j$ ($j=1,2$) and that
$$B_1B_2-2A_1C_2-2A_2C_1=t.\eqno (1.10)$$
It is easy to check that if $\tau\in SL_2({\bf Z})$, and
$\left(Q_1,Q_2\right)\in {\cal Q}_{d_1,d_2,t}$, then $\left(Q_1^{
\tau},Q_2^{\tau}\right)\in {\cal Q}_{d_1,d_2,t}$. Hence $SL_2({\bf Z}
)$
acts on ${\cal Q}_{d_1,d_2,t}$. Let us denote by $h\left(d_1,d_2,
t\right)$ the number
of $SL_2({\bf Z})$-equivalence classes of ${\cal Q}_{d_1,d_2,t}$. If $
t^2-d_1d_2\neq 0$,
then $h\left(d_1,d_2,t\right)$ is finite, it is proved in Appendix I of [M].

We now define the generalized class numbers in the
following way. If $d_1,d_2,t\in {\bf Z}$, $d_1\neq 0$, $d_2\neq 0$, $
t^2-d_1d_2\neq 0$ and $D_1,$$D_2$ are fundamental
discriminants with $D_i|d_i$ and $d_i/D_i\equiv 0,1$ (mod $4$) for $
i=1,2$, define
$$h_{D_1,D_2}\left(d_1,d_2,t\right):=\sum_{SL_2({\bf Z})\setminus
{\cal Q}_{d_1,d_2,t}}\omega_{D_1}\left(Q_1\right)\omega_{D_2}\left
(Q_2\right).\eqno (1.11)$$
If $\delta_1<0$, $\delta_{}{}_2<0$ are integers, let ${\cal R}_{\delta_{}{}_
1,\delta_2}$ be the subset of $\scr{Q}_{\delta_1}\times\scr{Q}_{\delta_
2}$ consisting of those pairs $\left(Q_1,Q_2\right)$ of quadratic forms satisfying that
$$Q_1=\lambda Q_2\hbox{\rm \ with some $\lambda\in {\bf Q}$}.$$
Note that ${\cal R}_{\delta_{}{}_1,\delta_2}$ is empty unless ${{
\delta_1}\over {\delta_2}}\in {\bf Q}^2.$ It is easy to check that if $
\tau\in SL_2({\bf Z})$, and
$\left(Q_1,Q_2\right)\in {\cal R}_{\delta_{}{}_1,\delta_2}$, then $\left
(Q_1^{\tau},Q_2^{\tau}\right)\in {\cal R}_{\delta_{}{}_1,\delta_2}$. Hence $
SL_2({\bf Z})$
acts on ${\cal R}_{\delta_{}{}_1,\delta_2}$. Let ${\cal R}_{\delta_
1,\delta_2}^{\ast}$ denote a complete set of representatives of the $
SL_2({\bf Z})$-equivalence classes of ${\cal R}_{\delta_{}{}_1,\delta_
2}$.

If $\delta_i<0$ are integers, $D_i$ are fundamental discriminants for $
i=1,2$ with $D_i|\delta_i$ and
$\delta_i/D_i\equiv 0,1$ (mod $4$), then define
$$E_{\delta_1,\delta_2,D_1,D_2}:=\sum_{\left(Q_1,Q_2\right)\in {\cal R}_{
\delta_1,\delta_2}^{\ast}}{{\omega_{D_1}\left(Q_1\right)\omega_{D_
2}\left(Q_2\right)}\over {\left|M\left(Q_1\right)\right|}}.\eqno
(1.12)$$
{\bf 1.5. Statement of the theorem.} Let $\beta >0$. We say that a
function $\chi$ satisfies Condition $A_{\beta}$ if $\chi$ is an even holomorphic function defined on the
strip $\left|\hbox{\rm Im}\,z\right|<\beta$ and the function
$$\left|\chi (z)\right|\left(1+\left|z\right|\right)^{\beta}$$
is bounded on this strip.

Let $F\left(\alpha ,\beta ,\gamma ;z\right)$ denote the Gauss hypergeometric function. If $
\chi\left(z\right)$ is a function for $z\ge 0$ and the following integral is absolutely convergent, introduce the notation
$$T_{\chi}\left(y\right):={1\over {2\pi}}\int_0^{\infty}\left|{{\Gamma\left
({1\over 4}+iz\right)\Gamma\left({3\over 4}+iz\right)}\over {\Gamma\left
(2iz\right)}}\right|^2F\left({1\over 4}-iz,{1\over 4}+iz,1,-y\right
)\chi\left(z\right)dz$$
for $y\ge 0$.

Let $\delta_{x,y}$ be Kronecker's symbol.

{\bf THEOREM 1.1.} {\it For} $i=1,2$ {\it let} $\delta_i<0$ {\it be integers. Let} $
D_i>0$
{\it be fundamental discriminants for} $i=1,2$ {\it with} $D_i|\delta_
i$ {\it and} $\delta_i/D_i\equiv 0,1$ (mod $4$). {\it There is an absolute constant} $
\beta >0$ {\it such that if} $\chi$ {\it is a function
satisfying condition} $A_{\beta}${\it , then the sum of}
$${{12}\over {\pi^2}}\delta_{1,D_1}\delta_{1,D_2}L^{\ast}\left(1,
\delta_1\right)L^{\ast}\left(1,\delta_2\right)\chi\left({i\over 4}\right
),$$
$$144\pi |\delta_1\delta_2|^{3/4}\sum_{j=1}^{\infty}\left(\hbox{\rm Shim}
F_j,\hbox{\rm Shim}F_j\right)_1b_j\left(D_1\right)\overline {b_j\left
({{\delta_1}\over {D_1}}\right)}b_j\left({{\delta_2}\over {D_2}}\right
)\overline {b_j\left(D_2\right)}\chi\left(r_j\right)$$
{\it and}
$$\int_{-\infty}^{\infty}{{L^{\ast}\left({1\over 2}-2i\rho ,D_1\right
)L^{\ast}\left({1\over 2}-2i\rho ,{{\delta_1}\over {D_1}}\right)L^{
\ast}\left({1\over 2}+2i\rho ,D_2\right)L^{\ast}\left({1\over 2}+
2i\rho ,{{\delta_2}\over {D_2}}\right)\chi\left(\rho\right)}\over {
\zeta\left(1+4i\rho\right)\zeta\left(1-4i\rho\right)}}d\rho$$
{\it equals}
$$E_{\delta_1,\delta_2,D_1,D_2}T_{\chi}\left(0\right)+\sum_{f\in
{\bf Z},f^2>\left|\delta_1\delta_2\right|}h_{D_1,D_2}\left(\delta_
1,\delta_2,f\right)T_{\chi}\left({{f^2}\over {\left|\delta_1\delta_
2\right|}}-1\right).$$
{\it Every summation and integral is absolutely convergent.}

{\bf REMARK 1.1.} In the special case when $D_1=D_2=1$, explicit
elementary expressions are given for the class numbers
$h_{1,1}\left(\delta_1,\delta_2,f\right)=h\left(\delta_1,\delta_2
,f\right)$ in [B6]. We expect that similar explicit
formulas can be proved for $h_{D_1,D_2}\left(\delta_1,\delta_2,f\right
)$ in the same
way also for general $D_i$.

{\bf REMARK 1.2.} The integral transform $\chi\rightarrow T_{\chi}\left
(y\right)$ is
well-known, it is a special case of the so-called Jacobi
transform, see e.g. [Ko]. Its inversion is also ecplicitly
known, therefore it is possible to state a formula also by
writing a general test function on the arithmetic side.

{\bf REMARK 1.3.} Observe that we have in fact a weighted
spectral sum of the product of four Fourier coefficients
of weight $1/2$, the weights being $\left(\hbox{\rm Shim}F_j,\hbox{\rm Shim}
F_j\right)_1$.

{\bf 1.6. Further notations.} In order to give a sketch of the
proof of the theorem in the next subsection, we have to
introduce the following notations. These notations will be
needed also later in the paper.

For $z,w\in\bbb H$ let
$$u(z,w)={{\left|z-w\right|^2}\over {4\hbox{\rm Im$z$Im$w$}}},\eqno
(1.13)$$
this is closely related to the hyperbolic distance $\rho (z,w)$ of $
z$
and $w$, namely we have $1+2u=\cosh\rho$.

If $m$ is a function on $[0,\infty )$, then for $z,w\in\bbb H$ write
$$m(z,w)=m\left(u(z,w)\right)\eqno (1.14)$$
by an abuse of notation. Conversely, if $m(z,w)$ is such a
function defined on $\bbb H\times\bbb H$ which depends only on
$u(z,w)$, then we can define a function $m$ on $[0,\infty )$ such
that (1.14) holds.

If $n,t$ are integers, $n>0$, let
$$\Gamma_{n,t}=\left\{\left(\matrix{a&b\cr
c&d\cr}
\right):\;a,b,c,d\in {\bf Z},\;ad-bc=n,\;a+d=t\right\}.$$
The group $SL_2({\bf Z})$ acts on this set by conjugation. If
$\gamma =$$\left(\matrix{a&b\cr
c&d\cr}
\right)\in\Gamma_{n,t}$, let $Q_{\gamma}(X,Y)=cX^2+(d-a)XY-bY^2$.
Then it is easy to see (see [B1], p 119) that this is a
one-to-one correspondence between $\Gamma_{n,t}$ and $\scr{Q}_{\delta}$ with $
\delta =t^2-4n$, and also
between the conjugacy classes of $\Gamma_{n,t}$ over $SL_2({\bf Z}
)$ and the
$SL_2({\bf Z})$-equivalence classes of $Q_{\delta}$. We remark that if
$\delta <0$, $\gamma\in\Gamma_{n,t}$, then $z_Q$$_{\gamma}$ is the unique fixed point of $
\gamma$ in $\bbb H$.

Let $n,t$ be integers, $n>0$, and for $\delta :=t^2-4n$ assume
$\delta\neq 0$. Let $D$ be a fundamental discriminant with $D|\delta$ and
$\delta /D\equiv 0,1$ (mod $4$). For a matrix $\gamma\in\Gamma_{n
,t}$ define
$$\omega_D\left(\gamma\right)=\omega_D\left(Q_{\gamma}(X,Y)\right
).$$
It is clear that if $\tau\in SL_2({\bf Z})$, then $\omega_D\left(
\tau^{-1}\gamma\tau\right)=\omega_D\left(\gamma\right)$.

If $E>0$, let ${\cal K}_E$ be the set of measurable functions $k$ on $
[0,\infty )$
satisfying that $k(u)\left(1+u\right)^E$ is bounded for $u\ge 0$.

Let $n,t$ be integers, $n>0$, and for $\delta :=t^2-4n$ assume
$\delta <0$. Let $D$ be a fundamental discriminant with $D|\delta$ and
$\delta /D\equiv 0,1$ (mod $4$). If $m\in {\cal K}_E$ for a large enough absolute
constant $E$, for $z\in\bbb H$ define
$$M_{t,n,D,m}(z):=\sum_{\gamma\in\Gamma_{n,t}}\omega_D\left(\gamma\right
)m\left(z,\gamma z\right).\eqno (1.15)$$
One can easily see using Lemma 2.1 below that this is a
bounded automorphic function on $\bbb H$.

Denote by $\Lambda_{\delta}$ a complete set of representatives of the
$SL_2({\bf Z})$-equivalence classes of $\scr{Q}_{\delta}$.

{\bf 1.7. Outline of the proof of Theorem 1.1.} We fix
integers $n_i$, $t_i$ for $i=1,2$ such that $\delta_i:=t_i^2-4n_i$. We take
two test functions $m_1,m_2\in {\cal K}_E$ for a large $E$ and consider the intagral
$$I:=\int_{{\cal F}_1}M_{t_1,n_1,D_1,m_1}(z)M_{t_2,n_2,D_2,m_2}(z
)d\mu_z.\eqno (1.16)$$
We compute $I$ in two different ways.

Firstly, just as in [B2, Lemma 2.2], using the definitions
of the functions $M_{t_1,n_1,D_1,m_1}(z)$ we give an elementary
expression for $I$ in Lemma 2.2 below involving the generalized
class numbers $h_{D_1,D_2}\left(\delta_1,\delta_2,f\right)$, where $
f$ runs over integers.

Secondly, we consider $I$ as the inner product of two
automorphic functions, and we compute this inner
product by the spectral theorem. To do so we have to consider integrals of the form
$$J_u:=\int_{{\cal F}_1}M_{t,n,D,m}(z)u(z)d\mu_z,\eqno (1.17)$$
where $u$ is a cusp form (or an Eisenstein series). We have
considered such integrals in our earlier papers, see [B1,
Lemma 2] and [B4, Lemma 3.2]. In the present case when
$\delta =t^2-4n$ is negative, the result is that
$$J_u=F_m\left(\lambda\right)\sum_{Q\in\Lambda_{\delta}}{{\omega_
D\left(Q\right)}\over {M_Q}}u\left(z_Q\right),\eqno (1.18)$$
where $F_m\left(\lambda\right)$ depends only on the given test function $
m$
and the Laplace-eigenvalue $\lambda$ of $u$ (considering $t$, $n$ and so $
\delta$ to be fixed). Now, by a Katok-Sarnak type formula the
summation $\sum_{Q\in\Lambda_{\delta}}{{\omega_D\left(Q\right)}\over {
M_Q}}u\left(z_Q\right)$ can be expressed
essentially as the product of {\it two\/} Fourier coefficients of
the cusp form $F$ of weight $1/2$ belonging to Kohnen's
subspace and satisfying that the Shimura lift of $F$ equals
$u$. The results of [I-L-T] and [B-M] will be important at this
step. When we compute $I$ by the spectral theorem, we have a
spectral sum of products of {\it two\/} integrals of the form
$J_u$. Therefore, finally we have a spectral sum of
products of {\it four\/} Fourier coefficients of weight $1/2$.

Choosing the test functions $m_i$ suitably we can get the theorem.

We will give the elementary expression for (1.16), and we
will express the integrals (1.17) of the form (1.18) in
Section 2. In Section 3 we compute the summations over
Heegner points occurring in (1.18). We complete the proof
of the theorem in Section 4.

{\bf 1.8. Discussion of the extension of Theorem 1.1 for other cases and statement of Theorem 1.2.}

The ideas sketched in Subsection 1.7 can be easily applied when we have more positive than negative Fourier coefficients. For example, when we have four positive Fourier coefficients, we consider the same integral (1.16) but in this case we have $\delta_1>0$, $\delta_2>0$. We can give an elementary expression for (1.16) extending the proof of [B2, Lemma 2.2]. The integral (1.17) is computed for this case in [B1, Lemma 2]. If there are three positive Fourier coefficients, we can still consider an integral of the form (1.16), but in this case we have to take $\delta_1>0$, $\delta_2<0$.

If there are more negative than positive Fourier coefficients, the same line of ideas can be still applied, but for this case we have to modify the definition of the function (1.15). We have to compute then the analogue of the integral (1.17).

This is done in Theorem 1.2 below, which is stated here and will be proved in Section 5. We note that for the proof of Theorem 1.2 the extension of the Katok-Sarnak formula for the case of two negative Fourier coefficients will be important. This extension was proved relatively recently in [D-I-T] and [I-L-T]. To state Theorem 1.2 we need the following notations.

If $Q(X,Y)=AX^2+BXY+CY^2$ is an element of $\scr{Q}_{\delta}$ with
some $\delta >0$, and $z_1$ and $z_2$ are the roots of $Az^2+Bz+C$ (if $
A=0$, one
root is $\infty$, otherwise these are real numbers), let $l_Q$ be
the noneuclidean line in $\bbb H$ connecting $z_1$ and $z_2$, let
$$C(Q)=\left\{\gamma\in PSL(2,{\bf Z}):\;\gamma z_1=z_1,\;\gamma
z_2=z_2\right\}.$$
If $A\neq 0$, this is an Euclidean semi-circle, and we orient
it counterclockwise for $A>0$, and clockwise for $A<0$. If
$A=0$ and $B>0$, then we orient the line $l_Q$ upwards, if
$A=0$ and $B<0$, then we orient it downwards. Finally let $C_Q=C(
Q)\setminus l_Q$, i.e we factorize by the action of $C(Q)$.

For $z,w\in\bbb H$ let
$$h(z,w):={{\left(z-\overline w\right)^2}\over {\left|z-\overline
w\right|^2}},\eqno (1.19)$$
see p 349 of [H] and also p 238 of [B5].

We now modify the definition (1.15) in the following way. Let $n,
t$ be integers, $n,t>0$, and for $\delta :=t^2-4n$ assume
$\delta >0$. Let $D$ be a fundamental discriminant with $D|\delta$ and
$\delta /D\equiv 0,1$ (mod $4$). If $m\in {\cal K}_E$ for a large enough absolute constant $
E$, for $z\in\bbb H$ define
$$N_{t,n,D,m}(z):=\sum_{\gamma\in\Gamma_{n,t}}\omega_D\left(\gamma\right
)m\left(z,\gamma z\right)h\left(\gamma z,z\right)\left({{j_{\gamma}
(z)}\over {\left|j_{\gamma}(z)\right|}}\right)^2.$$
For $\lambda <0$ consider the differential equation
$$f^{(2)}(\theta )={{\lambda}\over {\cos^2\theta}}f(\theta ),\qquad
\theta\in (-{{\pi}\over 2},{{\pi}\over 2}).\eqno (1.20)$$
Let $h_{\lambda}(\theta )$ be the unique odd solution of this equation with $
h_{\lambda}^{(1)}(0)=1$.

{\bf THEOREM 1.2.} {\it Let} $\delta >0$ {\it be an integer, let} $
D<0$ {\it be a fundamental discriminant with} $D|\delta$ {\it and} $
\delta /D\equiv 0,1$ (mod $4$).
{\it Let} $n,t$ {\it be positive integers such that} $t^2-4n=\delta$. {\it Let} $
u$ {\it be an even Hecke normalized Maass-Hecke cusp form for }
$SL_2({\bf Z})$ {\it with} $\Delta_0u=\lambda u$, $\lambda <0${\it , and let} $
u=\hbox{\rm Shim}F_j$ {\it for some }
$j\ge 1$. {\it Let} $m\in {\cal K}_E$ {\it with a large enough absolute constant} $
E${\it . Then we have}
$${1\over {\left(u,u\right)_1}}\int_{{\cal F}_1}N_{t,n,D,m}(z)u(z
)d\mu_z=\delta ^{3/4}\overline {b_j\left(D\right)}b_j\left({{\delta}\over
D}\right)F_{\delta ,n,m}\left(\lambda\right)$$
{\it with}
$$F_{\delta ,n,m}\left(\lambda\right):=48\sqrt {\pi}i\int_{-\pi /2}^{
\pi /2}m\left({{\delta}\over {4n\cos^2\theta}}\right){{\sqrt {1+{{
4n}\over {\delta}}}}\over {1+{{4n}\over {\delta}}\cos^2\theta}}h_{
\lambda}(\theta ){{\sin\theta d\theta}\over {\cos\theta}}.$$

{\bf 2. Inner product of automorphic functions}

In Subsection 2.1 we will give further notations needed
in Section 2 and we prove an upper bound, Lemma 2.1,
which will ensure that we will always have absolute
convergence in our calculations later. In Subsection 2.2 our main result is Lemma 2.2, which
gives an elementary expression for the integral $I$ defined
in (1.16) above. In Subsection 2.3 we express the
integrals $J_u$ given in (1.17) in the form (1.18).

{\bf 2.1. Notations and an upper bound. }

From now on, ${\cal F}_1$ will denote the closure of the standard fundamental
domain of the quotient $SL(2,{\bf Z})\setminus\bbb H$:
$${\cal F}_1:=\left\{z\in {\bf C}:\;\hbox{\rm Im$\,z>0$},\>-{1\over
2}\le\hbox{\rm Re$\,z$}\le{1\over 2},\>\left|z\right|\ge 1\right\}
.\eqno (2.1)$$
For $\phi\in [0,2\pi ]$, write
$$k_{\phi}=\left(\matrix{\cos\phi&\sin\phi\cr
-\sin\phi&\cos\phi\cr}
\right).\eqno (2.2)$$
These matrices form the stability group of $i$ in $SL_2({\bf R})$.

If $\gamma$ is an elliptic element of $PSL_2({\bf R})$, let
$$C\left(\gamma\right):=\left\{\tau\in SL_2({\bf Z}):\;\;\tau\gamma
=\gamma\tau\right\}.$$
It is well-known and easily proved that we have
$$C\left(\gamma\right)=\left\{\tau\in SL_2({\bf Z}):\;\;\tau z_{\gamma}
=z_{\gamma}\right\},\eqno (2.3)$$
where $z_{\gamma}$ is the unique fixed point of $\gamma$ in $\bbb
H$. It is also known that $C\left(\gamma\right)$ is always finite, it has an even number of elements, let $\left
|C_{\gamma}\right|=2M_{\gamma}$.

Note that the following lemma is a variant of Lemma 5.3 of [B2].

{\bf LEMMA 2.1.} {\it Let} $n,t$ {\it be integers,} $n>0$, $t^2-4
n<0$. {\it Let }
$z=x+iy\in {\cal F}_1$ {\it and} $X\ge 1${\it . Then for every} $
\epsilon >0$ {\it we have that}
$$\left|\left\{\gamma\in\Gamma_{n,t}:\;u\left(\gamma z,z\right)\le
X\right\}\right|\ll_{\epsilon ,t,n}X^{{1\over 2}+\epsilon}.$$
{\it Proof.\/} Let $\gamma\in\Gamma_{n,t}$, and write $\gamma =\left
(\matrix{a&b\cr
c&d\cr}
\right)$. First note that by [I], (1.9) and
(1.11) we have
$$4u\left(\gamma z,z\right)={{\left|cz^2+\left(d-a\right)z-b\right
|^2}\over {n\hbox{\rm Im$^2z$}}}.\eqno (2.4)$$
It is easy to compute that we have
$$\hbox{\rm Im}\left(cz^2+\left(d-a\right)z-b\right)=2cxy+\left(d
-a\right)y,$$
$$\hbox{\rm Re}\left(cz^2+\left(d-a\right)z-b\right)=c\left(x^2-y^
2\right)+\left(d-a\right)x-b.$$
Hence $u\left(\gamma z,z\right)\le X$ and (2.4) imply that
$$2cx+d-a\ll_n\sqrt X,\eqno (2.5)$$
$$c\left(x^2+y^2\right)+b\ll_n\sqrt Xy.\eqno (2.6)$$
We get from (2.5) that
$$d=-cx+O_{t,n}\left(\sqrt X\right),\;a=cx+O_{t,n}\left(\sqrt X\right
),$$
and from these relations and (2.6) we get
$$n=ad-bc=-c^2x^2+c^2\left(x^2+y^2\right)+O_{t,n}\left(\sqrt X\left
(\sqrt X+y\left|c\right|\right)\right).$$
This implies $c=O_{t,n}\left(\sqrt X\right)$, and so (2.5) gives
$d-a\ll_{t,n}\sqrt X.$ Then there are $O_{t,n}\left(\sqrt X\right
)$ possibilities for the
pair $\left(a,d\right)$. If $a$ and $d$ are given with $ad\neq n$, then
$bc=ad-n$ implies that there are $O_{\epsilon ,t,n}\left(X^{\epsilon}\right
)$ possibilities for the pair $\left(b,c\right)$.
Finally, if $ad=n$, then $a+d=t$ gives
$\left(a-d\right)^2=t^2-4n<0$, a contradiction. The lemma is proved.

{\bf 2.2. Inner product of two functions of type} $M_{t,n,D,m}(z)$ {\bf and pairs of quadratic forms.}

Let $-2<\hbox{\rm $\tau_1,$$\tau_2$}<2$ be real numbers and let $
m_1,m_2\in {\cal K}_E$ with
a large enough absolute constant $E>0$. For every $\Phi >1$ let us define
$$\!\!\!{\cal L}\left(\tau_1,\tau_2,\phi ,m_1,m_2\right):=\int\!\!\!
\int{{m_1\left(\left(4-\tau_1^2\right)r_1\left(r_1+1\right)\right
)m_2\left(\left(4-\tau_2^2\right)r_2\left(r_2+1\right)\right)dr_1
dr_2}\over {\sqrt {2\Phi\left(2r_1+1\right)\left(2r_2+1\right)-\Phi^
2-\left(2r_1+1\right)^2-\left(2r_2+1\right)^2+1}}},\eqno (2.7)$$
where we integrate over the set
$$\left\{\left(r_1,r_2\right)\in {\bf R}_{+}^2:\;2\Phi\left(2r_1+
1\right)\left(2r_2+1\right)-\Phi^2-\left(2r_1+1\right)^2-\left(2r_
2+1\right)^2+1\ge 0\right\}.\eqno (2.8)$$
Here ${\bf R}_{+}$ is the set of nonnegative real numbers.

{\bf LEMMA 2.2.} {\it For} $i=1,2$ {\it let} $n_i,t_i$ {\it be integers,} $
n_i>0${\it , and for }
$\delta_i:=t_i^2-4n_i$ {\it assume} $\delta_i<0${\it . Let} $D_i$ {\it be fundamental discriminants for} $
i=1,2$ {\it with} $D_i|\delta_i$ {\it and }
$\delta_i/D_i\equiv 0,1$  (mod $4$). {\it Let} $m_1,m_2\in {\cal K}_
E$ {\it with a large enough absolute constant} $E>0${\it .}

{\it Then using the notation (1.15) we have that}
$$\int_{{\cal F}_1}M_{t_1,n_1,D_1,m_1}(z)M_{t_2,n_2,D_2,m_2}(z)d\mu_
z\eqno (2.9)$$
{\it equals the sum of}
$$4\pi E_{\delta_1,\delta_2,D_1,D_2}\int_0^{\infty}m_1\left({{\left
|\delta_1\right|}\over {n_1}}r\left(1+r\right)\right)m_2\left({{\left
|\delta_2\right|}\over {n_2}}r\left(1+r\right)\right)dr$$
{\it and}
$$8\sum_{f\in {\bf Z},f^2>\left|\delta_1\delta_2\right|}h_{D_1,D_
2}\left(\delta_1,\delta_2,f\right){\cal L}\left({{t_1}\over {\sqrt {
n_1}}},{{t_2}\over {\sqrt {n_2}}},\left|{f\over {\sqrt {\left|\delta_
1\delta_2\right|}}}\right|,m_1,m_2\right).\eqno (2.10)$$
{\it The quantities} $E_{\delta_1,\delta_2,D_1,D_2}$ {\it and} $h_{
D_1,D_2}\left(\delta_1,\delta_2,f\right)$ {\it are defined in Subsection 1.4, the} $
{\cal L}${\it -function is defined in (2.7) and (2.8). The sum (2.10) is absolutely convergent.}

We postpone the proof of this lemma to the end of this
subsection. We first need three preliminary lemmas.

{\bf LEMMA 2.3.} {\it Use the notations and assumptions of Lemma 2.2. Write} $
G:=\Gamma_{n_1,t_1}\times\Gamma_{n_2,t_2},$ {\it and let} $G_0$ {\it be the set of those elements} $\left
(\gamma_1,\gamma_2\right)\in G$ {\it for}
{\it which the fixed point of} $\gamma_1$ {\it in} $\bbb H$ {\it coincides with the fixed point of} $
\gamma_2$ {\it in} $\bbb H${\it . If} $\left(\gamma_1,\gamma_2\right
),\left(\gamma_1^{\ast},\gamma_2^{\ast}\right)\in G${\it , we say that} $\left
(\gamma_1,\gamma_2\right)$ {\it and}
$\left(\gamma_1^{\ast},\gamma_2^{\ast}\right)$ {\it are} $SL_2({\bf Z}
)${\it -equivalent if there is an element} $\tau\in SL_2({\bf Z})$ {\it such that}
$\tau^{-1}\gamma_i\tau =\gamma_i^{\ast}$ {\it for} $i=1,2${\it . We denote by} $
G_0^{\ast}$ {\it a complete set of representatives of the} $SL_2(
{\bf Z})${\it -equivalence classes of} $G_0${\it , and by} $\left
(G\setminus G_0\right)^{\ast}$ {\it a complete set of representatives of the}
$SL_2({\bf Z})${\it -equivalence classes of} $G\setminus G_0${\it .}

{\it We have that}
$$\int_{{\cal F}_1}M_{t_1,n_1,D_1,m_1}(z)M_{t_2,n_2,D_2,m_2}(z)d\mu_
z\eqno (2.11)$$
{\it equals the sum of }
$$\sum_{\left(\gamma_1,\gamma_2\right)\in G_0^{\ast}}{{\omega_{D_
1}\left(\gamma_1\right)\omega_{D_2}\left(\gamma_2\right)}\over {M\left
(\gamma_1\right)}}\int_{\bbb H}m_1\left(z,\gamma_1z\right)m_2\left
(z,\gamma_2z\right)d\mu_z\eqno (2.12)$$
{\it and}
$$\sum_{\left(\gamma_1,\gamma_2\right)\in\left(G\setminus G_0\right
)^{\ast}}\omega_{D_1}\left(\gamma_1\right)\omega_{D_2}\left(\gamma_
2\right)\int_{\bbb H}m_1\left(z,\gamma_1z\right)m_2\left(z,\gamma_
2z\right)d\mu_z.\eqno (2.13)$$
{\it The integral (2.11) is absolutely convergent, and the integral and summation are absolutely convergent together in (2.12) and (2.13). }

{\it Proof.\/} Since $\delta_i<0$, any element $\gamma\in\Gamma_{
n_i,t_i}$ determines an
elliptic transformation of $\bbb H$, see Section 1.5 of [I].
Hence $\gamma$ has a unique fixed points in $\bbb H$. Assume that
$\gamma_1\in\Gamma_{n_1,t_1}$, $\gamma_2\in\Gamma_{n_2,t_2}$, $\tau
\in SL_2({\bf Z})$ and
$$\tau^{-1}\gamma_1\tau =\gamma_1,\quad\tau^{-1}\gamma_2\tau =\gamma_
2.\eqno (2.14)$$
It is clear by (2.3) that if $\left(\gamma_1,\gamma_2\right)\in G
\setminus G_0$, then (2.14) is true
if and only if $\tau =\pm\left(\matrix{1&0\cr
0&1\cr}
\right)$. If $\left(\gamma_1,\gamma_2\right)\in G_0$, then by (2.3)
we see that $C\left(\gamma_1\right)=C\left(\gamma_2\right)$, and (2.14) is true if and only
if $\tau\in C\left(\gamma_1\right)$. Recall that $C\left(\gamma_1\right
)$ is finite.

Therefore, if $\left(\gamma_1,\gamma_2\right)\in G\setminus G_0$, then the pairs
$$\left(\tau^{-1}\gamma_1\tau ,\tau^{-1}\gamma_2\tau\right)\eqno
(2.15)$$
$ $represent every element of the $SL_2({\bf Z})$-equivalence class of $\left
(\gamma_1,\gamma_2\right)$ exactly twice as $\tau$
runs over $SL_2({\bf Z})$. If $\left(\gamma_1,\gamma_2\right)\in
G_0$, then the pairs (2.15)
represent every element of the $SL_2({\bf Z})$-equivalence class
of $\left(\gamma_1,\gamma_2\right)$ exactly $\left|C\left(\gamma_
1\right)\right|$ times as $\tau$ runs over $SL_2({\bf Z})$.

By the definitions we see that (2.11) equals
$$\sum_{\gamma_1\in\Gamma_{t_1}}\sum_{\gamma_2\in\Gamma_{t_2}}\omega_
D\left(\gamma_1\right)\omega_D\left(\gamma_2\right)\int_{{\cal F}_
1}m_1\left(z,\gamma_1z\right)m_2\left(z,\gamma_2z\right)d\mu_z,$$
and Lemma 2.1 shows that the double summation and the
integration are absolutely convergent together. We partition $G$ into $
SL_2({\bf Z})$-equivalence classes. Since for
$\tau\in SL_2({\bf Z})$ we have that
$$\int_{{\cal F}_1}m_1\left(z,\tau^{-1}\gamma_1\tau z\right)m_2\left
(z,\tau^{-1}\gamma_2\tau z\right)d\mu_z=\int_{\tau {\cal F}_1}m_1\left
(z,\gamma_1z\right)m_2\left(z,\gamma_2z\right)d\mu_z,$$
our considerations above give the lemma.

{\bf LEMMA 2.4.} {\it Let} $\gamma =\left(\matrix{A&B\cr
C&D\cr}
\right)\in SL_2({\bf R})$ {\it be an elliptic element and let} $z\in\bbb H$ {\it be its fixed point. Let} $w\in\bbb H${\it . Then one has}
$$u\left(w,\gamma w\right)=4u\left(z,w\right)\left(u\left(z,w\right
)+1\right)C^2\hbox{\rm Im}^2z.$$
{\it Proof.\/} We use again the identity (as in (2.4))
$$u\left(w,\gamma w\right)={{\left|Cw^2+\left(D-A\right)w-B\right
|^2}\over {4\hbox{\rm $\hbox{\rm Im}^2w$}}}.$$
The roots of the quadratic polynomial $Cw^2+\left(D-A\right)w-B$
are $z$ and $\overline z$, hence
$$u\left(w,\gamma w\right)={{C^2\left|w-z\right|^2\left|w-\overline
z\right|^2}\over {4\hbox{\rm $\hbox{\rm Im}^2w$}}}.$$
One can check the identity
$$\left|w-\overline z\right|^2=\left|w-z\right|^2+4\hbox{\rm Im}z\hbox{\rm Im}
w.\eqno (2.16)$$
The lemma follows.

{\bf LEMMA 2.5.} {\it Let} $m_1,m_2\in {\cal K}_E$ {\it with a large enough absolute constant} $
E>0${\it . Let} $\gamma_1=\left(\matrix{a&b\cr
c&d\cr}
\right),$ $\gamma_2=\left(\matrix{A&B\cr
C&D\cr}
\right)$ {\it be elliptic elements of} $SL_2({\bf R})$. {\it Write} $
\tau_1=a+d${\it ,} $\tau_2=A+D${\it .}

{\it (i) Assume that} $\gamma_1$ {\it and} $\gamma_2$ {\it have different fixed points in} $\bbb
H${\it . Let}
$$F:=F\left(\gamma_1,\gamma_2\right)={{\left(d-a\right)\left(D-A\right
)+2bC+2Bc}\over {\sqrt {4-\left(d+a\right)^2}\sqrt {4-\left(D+A\right
)^2}}}.\eqno (2.17)$$
{\it Then we have} $|F|>1${\it , and (recalling (2.7) and (2.8)) we have that }
$$\int_{\bbb H}m_1\left(z,\gamma_1z\right)m_2\left(z,\gamma_2z\right
)d\mu_z=8{\cal L}\left(\tau_1,\tau_2,|F|,m_1,m_2\right).\eqno (2.
18)$$
{\it (ii) Assume that} $\gamma_1$ {\it and} $\gamma_2$ {\it have the same fixed point in }
$\bbb H${\it . Then we have that}
$$\int_{\bbb H}m_1\left(z,\gamma_1z\right)m_2\left(z,\gamma_2z\right
)d\mu_z=4\pi\int_0^{\infty}m_1\left(\left(4-\tau_1^2\right)r\left
(1+r\right)\right)m_2\left(\left(4-\tau_2^2\right)r\left(1+r\right
)\right)dr.\eqno (2.19)$$
{\it Proof.\/} First note that it is easy to check that
$F\left(\gamma_1,\gamma_2\right)=F\left(\tau^{-1}\gamma_1\tau ,\tau^{
-1}\gamma_2\tau\right)$ for $\tau\in SL_2({\bf R})$. The left-hand
sides of (2.18) and (2.19) also remain the same if we write $\tau^{
-1}\gamma_1\tau$ and
$\tau^{-1}\gamma_2\tau$ in place of $\gamma_1$ and $\gamma_2$, respectively. Therefore, it
is enough to prove the lemma for the pair $\left(\tau^{-1}\gamma_
1\tau ,\tau^{-1}\gamma_2\tau\right)$
with any $\tau\in SL_2({\bf R})$ instead of the pair $\left(\gamma_
1,\gamma_2\right)$.

Let $z_i$ be the fixed point of $\gamma_i$ in $\bbb H$ for $i=1,2$. We claim that there is a $
\sigma\in SL_2({\bf R})$ such that $\hbox{\rm $\hbox{\rm Im}\sigma$}
z_1=\hbox{\rm $\hbox{\rm Im}\sigma$}z_2$. Indeed, assume $\hbox{\rm Im}
z_1>\hbox{\rm Im}z_2$ and let $\sigma_d=\left(\matrix{0&-1\cr
1&d\cr}
\right)$ with some real $d$. Then
$$\hbox{\rm Im}\sigma_dz_i={{\hbox{\rm $\hbox{\rm Im}$}z_i}\over {\left
|z_i+d\right|^2}}.$$
If $d$ is large enough, then $\hbox{\rm Im}\sigma_dz_1>\hbox{\rm Im}
\sigma_dz_2$. If
$d=-\hbox{\rm $\hbox{\rm Re}z_2$}$, then
$$\hbox{\rm Im}\sigma_dz_2={{\hbox{\rm $\hbox{\rm Im}$}z_2}\over {\hbox{\rm Im}^
2z_2}}={1\over {\hbox{\rm Im}z_2}},$$
and so
$$\hbox{\rm Im}\sigma_dz_1\le{{\hbox{\rm $\hbox{\rm Im}$}z_1}\over {\hbox{\rm Im}^
2z_1}}={1\over {\hbox{\rm Im}z_1}}<{1\over {\hbox{\rm Im}z_2}}=\hbox{\rm Im}
\sigma_dz_2.$$
Therefore, by continuity there must be such a $d$ for
which $\hbox{\rm Im}\sigma_dz_1=\hbox{\rm Im}\sigma_dz_2$. Taking an appropriate upper triangular
element $\mu\in SL_2({\bf R})$ we can then achieve that
$$\hbox{\rm Im}\mu\sigma_dz_1=\hbox{\rm Im}\mu\sigma_dz_2=1,\quad\hbox{\rm Re}
\mu\sigma_dz_1=-\hbox{\rm Re}\mu\sigma_dz_2.$$
Hence replacing the pair $\left(\gamma_1,\gamma_2\right)$ with the pair $\left
(\tau^{-1}\gamma_1\tau ,\tau^{-1}\gamma_2\tau\right)$
for a suitable $\tau\in SL_2({\bf R})$ we can assume that
$$\hbox{\rm Im}z_1=\hbox{\rm Im}z_2=1,\quad\hbox{\rm Re}z_1=-\hbox{\rm Re}
z_2=X\eqno (2.20)$$
with some real $X$, where $z_i$ is the fixed point of $\gamma_i$ in $\bbb
H$ for $i=1,2$. In case (i) we have
$X\neq 0$, while in case (ii) we have $X=0$.

We assume (2.20) from now on.

The relation $\gamma_1z_1=z_1$ means
$$c\left(X+i\right)^2+\left(d-a\right)\left(X+i\right)-b=0,$$
which is equivalent to
$$2cX=a-d,\quad c\left(X^2+1\right)=-b.$$
Since $\tau_1=a+d$, we get $ad={{\tau_1^2}\over 4}-c^2X^2$. Then $
ad-bc=1$
implies $c^2=1-{{\tau_1^2}\over 4}$. We can compute every other entry
from $c$, we have $a={{\tau_1}\over 2}+cX$, $d={{\tau_1}\over 2}-
cX$, and finally we
have
$$\left(\matrix{a&b\cr
c&d\cr}
\right)=\left(\matrix{{{\tau_1}\over 2}+{{\epsilon_1X}\over 2}\sqrt {
4-\tau_1^2}&-{{\epsilon_1}\over 2}\sqrt {4-\tau_1^2}\left(X^2+1\right
)\cr
{{\epsilon_1}\over 2}\sqrt {4-\tau_1^2}&{{\tau_1}\over 2}-{{\epsilon_
1X}\over 2}\sqrt {4-\tau_1^2}\cr}
\right)$$
with some $\epsilon_1\in \{-1,1\}$. Similarly, we have
$$\left(\matrix{A&B\cr
C&D\cr}
\right)=\left(\matrix{{{\tau_2}\over 2}-{{\epsilon_2X}\over 2}\sqrt {
4-\tau_2^2}&-{{\epsilon_2}\over 2}\sqrt {4-\tau_2^2}\left(X^2+1\right
)\cr
{{\epsilon_2}\over 2}\sqrt {4-\tau_2^2}&{{\tau_2}\over 2}+{{\epsilon_
2X}\over 2}\sqrt {4-\tau_2^2}\cr}
\right)$$
with some $\epsilon_2\in \{-1,1\}$. For $z\in\bbb H$ we see by Lemma 2.4 that
$$u\left(z,\gamma_1z\right)=\left(4-\tau_1^2\right)u\left(z_1,z\right
)\left(u\left(z_1,z\right)+1\right),\eqno (2.21)$$
$$u\left(z,\gamma_2z\right)=\left(4-\tau_2^2\right)u\left(z_2,z\right
)\left(u\left(z_2,z\right)+1\right).\eqno (2.22)$$
Up to this point our reasoning is valid for both cases (i)
and (ii).

We now assume (i). Then by (2.17) we have that
$$|F|=2X^2+1.\eqno (2.23)$$
We get by (2.20) for $z=x+iy$ that
$$r_1:=u\left(z_1,z\right)={{\left(X-x\right)^2+\left(y-1\right)^
2}\over {4y}},\quad r_2:=u\left(z_2,z\right)={{\left(X+x\right)^2
+\left(y-1\right)^2}\over {4y}}.\eqno (2.24)$$
Then we have
$$r_2-r_1={{Xx}\over y}$$
and
$$r_2={{\left(X+{{r_2-r_1}\over X}y\right)^2+\left(y-1\right)^2}\over {
4y}},$$
which is the same as
$$0=2y\left(-r_2-r_1-1\right)+\left(1+\left({{r_2-r_1}\over X}\right
)^2\right)y^2+X^2+1.\eqno (2.25)$$
So if $X\neq 0$ and $r_1,r_2\ge 0$ are given, then there are real
numbers $x,y$ with $y>0$ satisfying (2.24) with $z=x+iy$ if
and only if
$$2X^2\left(2r_1r_2+r_1+r_2\right)-X^4-\left(r_2-r_1\right)^2\ge
0,\eqno (2.26)$$
and if this is true, then the pairs $\left(x,y\right)$ satisfying (2.24) are given by
$$y=y_1={{1+r_1+r_2+{1\over X}\sqrt {2X^2\left(2r_1r_2+r_1+r_2\right
)-X^4-\left(r_2-r_1\right)^2}}\over {1+\left({{r_2-r_1}\over X}\right
)^2}},\eqno (2.27)$$
$$y=y_2={{1+r_1+r_2-{1\over X}\sqrt {2X^2\left(2r_1r_2+r_1+r_2\right
)-X^4-\left(r_2-r_1\right)^2}}\over {1+\left({{r_2-r_1}\over X}\right
)^2}}\eqno (2.28)$$
and
$$x={{r_2-r_1}\over X}y.$$
By (2.24) we get
$${{dr_1}\over {dx}}={{x-X}\over {2y}},\;{{dr_2}\over {dx}}={{x+X}\over {
2y}},$$
$${{dr_1}\over {dy}}={1\over 4}-{{1+\left(X-x\right)^2}\over {4y^
2}},\;{{dr_2}\over {dy}}={1\over 4}-{{1+\left(X+x\right)^2}\over {
4y^2}}.$$
Hence we can compute that
$$\left|{{dr_1dr_2}\over {dxdy}}\right|=\left|{{X\left(1+X^2-x^2-
y^2\right)}\over {4y^3}}\right|.\eqno (2.29)$$
We have by (2.24) that
$${{1+X^2}\over y}-1-r_1-r_2={{1+X^2-x^2-y^2}\over {2y}}.\eqno (2
.30)$$
We have by (2.25) that the product of the two roots of
that quadratic polynomial in $y$ is $y_1y_2={{1+X^2}\over {1+\left
({{r_2-r_1}\over X}\right)^2}}$.
Alternatively, we can see it directly from (2.27) and
(2.28). Hence for $i=1,2$ we have
$${{1+X^2}\over {y_i}}-1-r_1-r_2=\left(1+\left({{r_2-r_1}\over X}\right
)^2\right)y_{3-i}-1-r_1-r_2,$$
hence (2.27) and (2.28) give that
$$\left|{{1+X^2}\over {y_i}}-1-r_1-r_2\right|=\left|{1\over X}\sqrt {
2X^2\left(2r_1r_2+r_1+r_2\right)-X^4-\left(r_2-r_1\right)^2}\right
|\eqno (2.31)$$
for $i=1,2.$ Then (2.29), (2.30) and (2.31) show that
$$\left|{{dxdy}\over {y^2}}\right|=\left|{{2dr_1dr_2}\over {\sqrt {
2X^2\left(2r_1r_2+r_1+r_2\right)-X^4-\left(r_2-r_1\right)^2}}}\right
|.$$
Substituting $\left(r_1,r_2\right)$ in place of $\left(x,y\right)$ by (2.24), we get by
(2.26), (2.27), (2.28), (2.21) and (2.22), that the left-hand side of (2.18) equals
$$4\int\!\!\!\int{{m_1\left(\left(4-\tau_1^2\right)r_1\left(r_1+1\right
)\right)m_2\left(\left(4-\tau_2^2\right)r_2\left(r_2+1\right)\right
)}\over {\sqrt {2X^2\left(2r_1r_2+r_1+r_2\right)-X^4-\left(r_2-r_
1\right)^2}}}dr_1dr_2,$$
where we integrate over the set
$$\left\{\left(r_1,r_2\right)\in {\bf R}_{+}^2:\;\;2X^2\left(2r_1
r_2+r_1+r_2\right)-X^4-\left(r_2-r_1\right)^2\ge 0\right\}.$$
Taking into account (2.23) we obtain part (i) of the lemma.

Now consider case (ii). We then take geodesic polar
coordinates around $i$: for every $z\in\bbb H$ we can uniquely write
$${{z-i}\over {z+i}}=\tanh\left({R\over 2}\right)e^{i\phi}$$
with $R>0$ and $0\le\phi <2\pi$. It is known and easily computed
that the invariant measure is expressed in these new
coordinates as $d\mu_z=\sinh RdRd\phi .$ It follows from (2.16) that
$${1\over {\tanh^2\left({R\over 2}\right)}}=1+{1\over {u(z,i)}},$$
hence $u(z,i)=\sinh^2\left({R\over 2}\right)$. In case (ii) we have $
z_1=z_2=i$,
so using (2.21), (2.22) we get that the left-hand side of (2.19) equals
$$\int_0^{\infty}\int_0^{2\pi}m_1\left(\left(4-\tau_1^2\right)r\left
(R\right)\left(1+r\left(R\right)\right)\right)m_2\left(\left(4-\tau_
2^2\right)r\left(R\right)\left(1+r\left(R\right)\right)\right)\sinh
RdRd\phi ,$$
where $r=r\left(R\right):=\sinh^2\left({R\over 2}\right)$. Substituting $
r$ in place of $R$
we get ${{dr}\over {dR}}={{\sinh R}\over 2}$. The lemma is proved.

{\it Proof of Lemma 2.2.\/} We apply Lemma 2.3, Lemma 2.5 and the bijection $
\gamma\rightarrow Q_{\gamma}$ between $\Gamma_{n_i,t_i}$ and $\scr{
Q}_{\delta_i}$
described in Subsection 1.6. Note that if
$\left(\matrix{a&b\cr
c&d\cr}
\right)\in\Gamma_{n_1,t_1},$ $\left(\matrix{A&B\cr
C&D\cr}
\right)\in\Gamma_{n_2,t_2}$, then we apply Lemma 2.5
for $\gamma_1=\left(\matrix{a/\sqrt {n_1}&b/\sqrt {n_1}\cr
c/\sqrt {n_1}&d/\sqrt {n_1}\cr}
\right)$, $\gamma_2=\left(\matrix{A/\sqrt {n_2}&B/\sqrt {n_2}\cr
C/\sqrt {n_2}&D/\sqrt {n_2}\cr}
\right)$. Recall the formulas (1.9)-(1.12). The lemma is proved.

{\bf 2.3. Spectral coefficients of functions of type} $M_{t,n,D,m}
(z)${\bf .}

As in [B3], for $\lambda <0$ let $g_{\lambda}(r)$ ($r\in [0,\infty
)$) be the unique solution of
$$g^{(2)}(r)+{{\cosh r}\over {\sinh r}}g^{(1)}(r)=\lambda g(r)\eqno
(2.32)$$
with $g_{\lambda}(0)=1.$ Writing $\lambda =-{1\over 4}-\tau^2$ with a complex $
\tau$ one can check the explicit formula
$$g_{\lambda}(r)=F\left({1\over 2}+i\tau ,{1\over 2}-i\tau ,1;-\sinh^
2{r\over 2}\right)\eqno (2.33)$$
for $r\ge 0$. Indeed, writing $g(r)=F\left(\sinh^2{r\over 2}\right
)$ with a
function $F\left(u\right)$ (as in [I], (1.20), (1.21)) defined for $
u\in [0,\infty )$ the differential
equation (2.32) becomes
$u\left(1+u\right)F^{(2)}(u)+\left(1+2u\right)F^{(1)}(u)=\lambda
F(u)$. This equation is
discussed on [I], pp 26-27 and it is shown there that
the only solution with $F(0)=1$ is $F\left({1\over 2}+i\tau ,{1\over
2}-i\tau ,1;-u\right)$.
Note that there is a misprint there in the displayed
formula between (1.43) nad (1.44), $-u$ should be there in
place of $u$. Let us define $g_0(r)=1$ for every $r\ge 0$. Then
(2.33) is true for every $\lambda\le 0$ and $r\ge 0$.

Every step of the proof of the next lemma can be found in the papers [B4], [B3], but for the sake of completeness we give the full proof here.

{\bf LEMMA 2.6.} {\it Let} $n,t$ {\it be integers,} $n>0${\it , write} $
\delta =t^2-4n$ {\it and assume} $\delta <0${\it . Let} $D$ {\it be a fundamental discriminant with} $
D|\delta$ {\it and }
$\delta /D\equiv 0,1$ {\it (mod} $4${\it ). Let} $m\in {\cal K}_E$ {\it with a large enough absolute constant} $
E${\it . Let} $u$ {\it be a Maass form of weight} $0$
{\it on} $\bbb H$ {\it and assume that} $\int_{{\cal F}_1}|u(z)|d
\mu_z<\infty${\it . Let} $\Delta_0u=\lambda u$ {\it with} $\lambda
\le 0${\it . Then we have }
$$\int_{{\cal F}_1}M_{t,n,D,m}(z)u(z)d\mu_z=\left(\sum_{Q\in\Lambda_{
\delta}}{{2\pi\omega_D\left(Q\right)}\over {M_Q}}u\left(z_Q\right
)\right)\int_0^{\infty}m\left({{\left|\delta\right|}\over {4n}}\sinh^
2r\right)g_{\lambda}(r)\sinh rdr.\eqno (2.34)$$
{\it If} $D>0$ {\it and} $u\left(z\right)=-u\left(-\overline z\right
)$ {\it for every} $z\in\bbb H${\it , then the left-hand side of (2.34) is} $
0${\it .}

{\it Proof.\/} We first prove (2.34). We see by (1.15) that the
left-hand side of (2.34) equals
$$\sum_{\gamma\in\Gamma_{n,t}}\omega_D\left(\gamma\right)\int_{{\cal F}_
1}m\left(z,\gamma z\right)u(z)d\mu_z,$$
and Lemma 2.1 and $\int_{{\cal F}_1}|u(z)|d\mu_z<\infty$ show that the summation and the
integration are absolutely convergent together.

We partition $\Gamma_{n,t}$ into conjugacy classes over $SL_2({\bf Z}
)$, for $\gamma\in\Gamma_{n,t}$ let
$$[\gamma ]=\left\{\tau^{-1}\gamma\tau :\,\tau\in SL_2({\bf Z})\right
\}.$$
If, for any $\gamma\in\Gamma_{n,t}$, we write
$$T_{\gamma}=\sum_{\delta\in [\gamma ]}\int_{{\cal F}_1}m\left(z,
\delta z\right)u(z)d\mu_z,$$
then we have
$$T_{\gamma}=\int_{C(\gamma )\setminus\bbb H}m\left(z,\gamma z\right
)u(z)d\mu_z.$$
Choose $h\in SL_2({\bf R})$ such that $h(i)=z_{\gamma}$, where $z_{
\gamma}$ is the
fixed point of $\gamma$ in $\bbb H$. Then recalling (2.2) there is a
$\phi_{\gamma}\in [0,\pi ]$ such that
$$h^{-1}\gamma hz=k_{\phi_{\gamma}}z\eqno (2.35)$$
for every $z\in\bbb H$. We get
$$T_{\gamma}={1\over {M_{\gamma}}}\int_{\bbb H}m\left(z,k_{\phi_{
\gamma}}z\right)u(hz)d\mu_z.$$
We use the substitution $z=k_{\phi}e^{-r}i$,
i.e. we use geodesic polar coordinates (see [I], Section
1.3), where $r\in\left(0,\infty\right)$, $\phi\in\left(0,\pi\right
)$. We have
$d\mu_z=\left(2\sinh r\right)drd\phi$, so using (1.13), (1.14) and also that $
k_{\phi_{\gamma}}$ and $k_{\phi}$
commute we get
$$T_{\gamma}={1\over {M_{\gamma}}}\int_0^{\infty}m\left(\left(\sin^
2\phi_{\gamma}\right)\sinh^2r\right)\left(\int_0^{\pi}u\left(h\left
(k_{\phi}e^{-r}i\right)\right)d\phi\right)\left(2\sinh r\right)dr
.$$
Let us define
$$G\left(z\right):=\int_0^{\pi}u\left(h\left(k_{\phi}z\right)\right
)d\phi$$
for $z\in\bbb H$. One obtains $G\left(z\right)$ by averaging the function
$u\left(hz\right)$ over the stability group of $i$ in $SL_2({\bf R}
)$, so $G(z)$ is radial at i, i.e. it depends only on the
 noneuclidean distance of $z$ and $i$ (see [I], Lemma 1.10). On the other hand,
 since $u$ is an eigenfunction of $\Delta_0$ with eigenvalue $\lambda$, so
is $G\left(z\right)$, because $\Delta_0$ commutes with the group action. A
radial (at $i$) eigenfunction of $\Delta_0$ with eigenvalue $\lambda$ is
determined up to a constant factor ([I], Lemma 1.12), so
using the form of the Laplace operator in geodesic polar
coordinates (see [I], (1.20)) and recalling (2.32) we get that
$$G\left(e^{-r}i\right)=\pi u\left(z_{\gamma}\right)g_{\lambda}\left
(r\right),$$
since $h(i)=z_{\gamma}$. We obtain
$$T_{\gamma}={{2\pi}\over {M_{\gamma}}}u\left(z_{\gamma}\right)\int_
0^{\infty}m\left(\left(\sin^2\phi_{\gamma}\right)\sinh^2r\right)g_{
\lambda}(r)\sinh rdr.$$
It follows from (2.35) and $\gamma\in\Gamma_{n,t}$ that $\left|2\cos
\phi_{\gamma}\right|=\left|{t\over {\sqrt n}}\right|$, so $\sin^2
\phi_{\gamma}={{\left|\delta\right|}\over {4n}}$. By the
remarks in Subsection 1.6 on the correspondence between $\Gamma_{
n,t}$
and $Q_{\delta}$ we obtain (2.34).

We now show the last statement of the lemma. It is not
hard to check that if $\gamma =\left(\matrix{a&b\cr
c&d\cr}
\right)\in\Gamma_{n,t}$, then we have
$\gamma^{\ast}:=\left(\matrix{a&-b\cr
-c&d\cr}
\right)\in\Gamma_{n,t}$ and
$$m\left(z,\gamma z\right)=m\left(-\overline z,\gamma^{\ast}\left
(-\overline z\right)\right)$$
for every $z\in\bbb H$. Since $D>0$, we have $\left({D\over {-1}}\right
)=1$, see
[D], p 41. This gives $\omega_D\left(\gamma\right)=\omega_D\left(
\gamma^{\ast}\right)$ by the definitions.
Hence we have $M_{t,n,D,m}(z)=M_{t,n,D,m}(-\overline z)$ for every
$z\in\bbb H$. Taking into account (2.1) the lemma follows.

{\bf 3. Shimura lifts, Zagier} $L-${\bf functions, Heegner points}

Our main result in this section is Lemma 3.5, where we
express the sum of Maass forms of weight $0$ over
Heegner points of a given discriminant. First we
analyze the Shimura lift in detail in Subsection 3.1, this
will be needed to handle the case of cusp forms. Then
we prove an elementary identity in Subsection 3.2, which
will be needed for the case of Eisenstein series.

{\bf 3.1. On Shimura lifts.} Let $\hbox{\rm $F\in V$$^{+}$}$ be a Maass cusp form of weight ${
1\over 2}$ satisfying
$\Delta_{1/2}F=s(s-1)F$ with some $s={1\over 2}+it$ and having the Fourier expansion
$$F(z)=\sum_{m\neq 0,\,m\equiv 0,1(4)}b_F(m){\rm W}_{{1\over 4}{\rm s}
{\rm g}{\rm n}(m),it}\left(4\pi\left|m\right|y\right){\rm e}\left
(mx\right)$$
for $z=x+iy\in\bbb H$.

Let $d$ be a fundamental discriminant, then we define the $d$th Shimura lift of $
F$ by
$$Sh_dF\left(z\right)=\sum_{k\neq 0}a_{Sh_dF}\left(k\right)W_{0,2
it}(4\pi\left|k\right|y)e(kx),\eqno (3.1)$$
where
$$a_{Sh_dF}\left(k\right):=\sum_{PQ=k,P>0}{{\left|Q\right|^{{1\over
2}}}\over P}\left({d\over P}\right)b_F\left(dQ^2\right).\eqno (3.
2)$$
Then it is known that $Sh_dF$ is an even weight $0$ cusp form for
the group $SL_2({\bf Z})$, see the proof of Proposition 6
(especially the lines below formula (10.6)) in
[D-I-T]; note that for $d>0$ it is also proved in Theorem 1 of
[B1].

If $0\neq F\in V$$^{+}$, then there is a fundamental
discriminant $d$ such that $Sh_dF$ is nonzero. Indeed, if
$Sh_dF$$=0$ for every fundamental discriminant $d$, then the
right-hand side of (3.2) is $0$ for every integer $k\ge 1$ and
for every fundamental discriminant $d$. Applying Mobius
inversion for a given $d$ we see that $b_F\left(dQ^2\right)=0$ for
every fundamental discriminant $d$ and for every integer
$Q$. It is not hard to see that every integer $n\equiv 0,1(4)$ can
be written in this form, so $b_F\left(n\right)=0$ for every such $
n$,
hence for every $n$, i.e. $F=0$, a contradiction.

Introduce the weight $0$ Hecke operators for every positive integer $
n$:
$$\left(H_nF\right)(z)={1\over {\sqrt n}}\sum_{ad=n,\;b\,{\rm m}{\rm o}
{\rm d}\,d}F\left({{az+b}\over d}\right),$$
where $a$ and $d$ run over positive integers. The Hecke operators $T_{p^2}$ of weight $1/2$ are defined in [K-S], p 199, see also our Subsection 1.2.

{\bf LEMMA 3.1.} {\it Let} $\hbox{\rm $F\in V$$^{+}$}$ {\it be a Maass cusp form of weight }
${1\over 2}$ {\it with} $\Delta_{1/2}F=\left(-{1\over 4}-t^2\right
)F$. {\it Let} $d$ {\it be a fundamental discriminant. Then\/} we have
$$\Delta_0\left(Sh_dF\right)=\left(-{1\over 4}-4t^2\right)Sh_dF,\eqno
(3.3)$$
{\it and for any prime} $p>2$ {\it we have that}
$$Sh_d\left(T_{p^2}F\right)=H_p\left(Sh_dF\right).\eqno (3.4)$$
{\it Proof.\/} Formula (3.3) follows at once from (3.1).

We prove (3.4) by showing that the Fourier
coefficients of both sides are the same. This can be done,
since we know the action of the operators on Fourier
coefficients: for $H_p$ see (1.1) of [K-S]; for $Sh_d$ see (3.2)
above; for $T_{p^2}$ see (1.3) of [K-S].  Using these formulas and that
$dQ^2\equiv 0,1$(4) is always true, we see that for
any integer $k\neq 0$ the $k$th Fourier coefficient of the
left-hand side of (3.4) is
$$\sum_{PQ=k,P>0}{{\left|Q\right|^{{1\over 2}}}\over P}\left({d\over
P}\right)\left(pb_F\left(dQ^2p^2\right)+p^{-1/2}\left({{dQ^2}\over
p}\right)b_F\left(dQ^2\right)+p^{-1}b_F\left({{dQ^2}\over {p^2}}\right
)\right),\eqno (3.5)$$
and the $k$th Fourier coefficient of the right-hand side of
(3.4) is
$$p^{1/2}\sum_{PQ=kp,P>0}{{\left|Q\right|^{{1\over 2}}}\over P}\left
({d\over P}\right)b_F\left(dQ^2\right)+p^{-1/2}\sum_{PQ=k/p,P>0}{{\left
|Q\right|^{{1\over 2}}}\over P}\left({d\over P}\right)b_F\left(dQ^
2\right).\eqno (3.6)$$
As in [K-S], we mean that $b_F\left(t\right)=0$ if $t$ is not an integer.

Now, the first term of (3.5) gives the $p$ divides $Q$ part of
the first term of (3.6), and the second term of (3.5) gives
the $p$ does not divide $Q$ part of
the first term of (3.6). Finally, the third term of (3.5)
equals the second term of (3.6). To see this we note that
$b_F\left({{dQ^2}\over {p^2}}\right)\neq 0$ implies that $p$ divides $
Q$. Indeed, since $p$ is
odd, ${{dQ^2}\over {p^2}}$ cannot be an integer if $p$ does not divide
$Q$, because the fundamental discriminant $d$ is not divisible
by $p^2$. So finally (3.5) equals (3.6), the lemma is proved.

Let $j\ge 1$ be given. Take a fundamental discriminant $d$ such that
$Sh_dF_j\neq 0$. By Lemma 3.1 we then get that $Sh_dF_j$ is a
weight $0$ Maass-Hecke cusp form for $SL(2,{\bf Z})$ whose $p$th
Hecke-eigenvalue is the $T_{p^2}$-eigenvalue of $F_j$ for every
prime $p>2$. By the Strong Multiplicity One Theorem it
follows that the first Fourier coefficient of $Sh_dF_j$ is
nonzero, i.e. (using (3.2)) we get that $b_j\left(d\right)\neq 0$. Let us
define
$$\hbox{\rm Shim}F_j\left(z\right):={1\over {b_j\left(d\right)}}S
h_dF_j\left(z\right).$$
Using again Lemma 3.1 and the Strong Multiplicity One Theorem we see
that this is well-defined (i.e. we get the same function
using any fundamental discriminant $d$ such that
$Sh_dF_j\neq 0$). Note that $\hbox{\rm Shim}F_j$ is an even Hecke normalized
Maass-Hecke cusp form of weight $0$ for $SL_2({\bf Z})$.

{\bf LEMMA 3.2.} {\it (i) The map} $j\rightarrow\hbox{\rm Shim}F_
j$ {\it gives a bijection between the positive integers and the even Hecke normalized Maass-Hecke cusp forms of weight} $
0$ {\it for} $SL_2({\bf Z})${\it .}

{\it (ii) If} $j\ge 1$ {\it is an integer,} $d$ {\it is a fundamental discriminant and for some} $
F\in V$$^{+}$ {\it we have} $Sh_dF=c\hbox{\rm Shim}F_j$ {\it with some} $
c\neq 0${\it , then} $F$ {\it is a constant multiple of} $F_j${\it . }

{\it Proof.\/} We first prove (i). We have seen above that this map is well-defined.
The injectivity of the map follows from our Lemma 3.1
and from Theorem 1.2 of [B-M].

To see the surjectivity first claim that if a cusp form $0\neq F\in
V^{+}$ is a common eigenfunction of $\Delta_{{1\over 2}}$ and the Hecke operators
$T_{p^2}$ for all but finitely many primes $p$, then $F$ is a
constant multiple of one of the basis elements $F_j$. Indeed, if
$$F=\sum_{j=1}^{\infty}c_jF_j,\eqno (3.7)$$
then $c_j\neq 0$ for a given $j$ implies that if $F$ is the
eigenfunction of a given $T_{p^2}$, then the $T_{p^2}$-eigenvalue of
$F_j$ is the same as that of $F$. The injectivity of the map
$j\rightarrow\hbox{\rm Shim}F_j$ and the Strong Multiplicity One Theorem then
implies that there may be only one $j$ for which $c_j\neq 0$
in (3.7).

Now, the surjectivity of the map $j\rightarrow\hbox{\rm Shim}F_j$ follows
from this claim applying again Theorem 1.2 of
[B-M]. Part (i) is proved.

Part (ii) follows easily from our Lemma 3.1 and from
Theorem 1.2 of [B-M]. The lemma is proved.

{\bf 3.2. An elementary identity.} Let $D$ be a fundamental discriminant and let $
\delta\neq 0$ be an
integer such that $D|\delta$ and $\delta /D\equiv 0,1$ {\it (mod} $
4${\it ).\/} For every
positive integer $q$ define
$$\rho_q\left(D,\delta\right):=\sum_{r\,{\rm m}{\rm o}{\rm d}\,2q
,\;r^2\equiv\delta\left(4q\right)}\omega_D\left(qX^2+rXY+{{r^2-\delta}\over {
4q}}Y^2\right).\eqno (3.8)$$
{\bf LEMMA 3.3.} {\it If} $D$ {\it is a fundamental discriminant,} $
\delta\neq 0$ {\it is an integer such that} $D|\delta$ {\it and} $
\delta /D\equiv 0,1$ {\it (\/}mod $4$){\it , then for every integer} $
q$ {\it we have that}
$$\sum_{q_1q_2=q}\mu\left(q_2\right)\left({D\over {q_2}}\right)\rho_{
q_1}\left(D,\delta\right)=\sum_{q_1q_2=q}\mu\left(q_2\right)\rho_{
q_1}\left(1,{{\delta}\over D}\right).\eqno (3.9)$$
{\it Proof.\/} One can give a function $f:{\bf Z}^2\rightarrow {\bf C}$ such that for every positive integer $
q$ we have that
$$\rho_q\left(D,\delta\right)=\sum_{d|q}\left({D\over {q/d}}\right
)f\left({{\delta}\over D},d\right).\eqno (3.10)$$
This follows from Theorem A of [B1], which is in
fact a reformulation of [K1, Proposition 5]. Indeed, we
apply Theorem A of [B1] with $T=0$, $c=q$, $s=\delta$,
$\hat {c}=4d$, noting that we have a nonzero term in the
second summation in Theorem A of [B1] only in
case $4|\hat {c}$.

Then by (3.10) we get that the left-hand side of (3.9) equals
$$\sum_{q_1q_2=q}\mu\left(q_2\right)\left({D\over {q_2}}\right)\sum_{
d|q_1}\left({D\over {q_1/d}}\right)f\left({{\delta}\over D},d\right
),$$
and writing $e:=q_1/d$ and $E:=q/d$ this equals
$$\sum_{dE=q}f\left({{\delta}\over D},d\right)\left({D\over E}\right
)\sum_{q_2e=E}\mu\left(q_2\right).$$
The inner sum is $0$ unless $E=1$, hence we proved that
$$\sum_{q_1q_2=q}\mu\left(q_2\right)\left({D\over {q_2}}\right)\rho_{
q_1}\left(D,\delta\right)=f\left({{\delta}\over D},q\right)\eqno
(3.11)$$
for every $D$, $\delta$ and $q$ satisfying the conditions of the
lemma. Applying (3.11) writing $1$ in place of $D$ and ${{\delta}\over
D}$ in
place of $\delta$ we obtain
$$\sum_{q_1q_2=q}\mu\left(q_2\right)\rho_{q_1}\left(1,{{\delta}\over
D}\right)=f\left({{\delta}\over D},q\right).\eqno (3.12)$$
The lemma follows from (3.11) and (3.12).

{\bf 3.3. Summation over Heegner points.}

Let $E\left(z,s\right)$ be the Eisenstein series for $PSL_2({\bf Z}
)$, see [I], Chapter 3.

{\bf LEMMA 3.4.} {\it Let} $D>0$ {\it be a fundamental discriminant and let} $
\delta <0$ {\it be an integer. Assume that} $D|\delta$ {\it and} $
\delta /D\equiv 0,1$ {\it (mod} $4${\it ). If} $\hbox{\rm Re$s$}>
1${\it , then}
$${1\over 2}\sum_{Q\in\Lambda_{\delta}}{{\omega_D\left(Q\right)}\over {
M_Q}}E\left(z_Q,s\right)=\left({{|\delta |}\over 4}\right)^{s/2}\sum_{
q=1}^{\infty}{{\rho_q\left(D,\delta\right)}\over {q^s}}.\eqno (3.
13)$$
{\it Proof.\/} This follows from Proposition 3.6 of [I-L-T]. We apply that proposition with $
k=0$,
$m=0$, $N=1$ (i.e we take there the group $\Gamma =SL(2,{\bf Z})$).
Our $D$ is denoted by $d$ there, and our $\delta$ is denoted by $
D$
there. The left-hand side of (3.13) equals the left-hand
side of the displayed equation in Proposition 3.6 of [I-L-T], since only the equivalence classes of
positive definite quadratic forms are considered there
(it can be seen a few lines above [I-L-T, Definition 1.2]),
while we consider both positive definite and negative definite
forms. We use also that $D>0$ implies $\left({D\over {-1}}\right)
=1$, see [D], p
41. The right-hand sides are also the same, taking
into account that in (3.4) of [I-L-T] the
variable $b$ runs modulo $c$, and not modulo $c/2$. The
lemma is proved.

{\bf LEMMA 3.5.} {\it Let} $D>0$ {\it be a fundamental discriminant and let} $
\delta <0$ {\it be an integer. Assume that} $D|\delta$ {\it and} $
\delta /D\equiv 0,1$ {\it (mod} $4${\it ). }

{\it (i) If} $\hbox{\rm Re$s$}={1\over 2}${\it , then}
$$\sum_{Q\in\Lambda_{\delta}}{{\omega_D\left(Q\right)}\over {M_Q}}
E\left(z_Q,s\right)=2\left({{|\delta |}\over 4}\right)^{s/2}{{L\left
(s,D\right)L\left(s,{{\delta}\over D}\right)}\over {\zeta\left(2s\right
)}}.$$
{\it (ii) We have that }
$$\sum_{Q\in\Lambda_{\delta}}{{\omega_D\left(Q\right)}\over {M_Q}}
=\delta_{1,D}{{\pi}\over {3\zeta\left(2\right)}}|\delta |^{1/2}L\left
(1,\delta\right),$$
{\it where} $\delta_{1,D}$ {\it is the Kronecker symbol. }

{\it (iii) If} $u$ {\it is an even Hecke normalized Maass-Hecke cusp form for} $SL_2({\bf Z})$ {\it and} $u=\hbox{\rm Shim}F_j$ {\it for some} $
j\ge 1${\it , then}
$${1\over {\left(u,u\right)_1}}\sum_{Q\in\Lambda_{\delta}}{{\omega_
D\left(Q\right)}\over {M_Q}}u\left(z_Q\right)=12|\delta |^{3/4}\overline {
b_j\left(D\right)}b_j\left({{\delta}\over D}\right).\eqno (3.14)$$
{\it Proof.\/} Lemma 3.3 gives that we have
$$\sum_{q=1}^{\infty}{{\rho_q\left(D,\delta\right)}\over {q^s}}={{
L\left(s,\chi_D\right)}\over {\zeta\left(s\right)}}\sum_{q=1}^{\infty}{{
\rho_q\left(1,{{\delta}\over D}\right)}\over {q^s}}.$$
for $\hbox{\rm Re$s$}>1$. Hence from (3.8) and (1.4) we get
$$\sum_{q=1}^{\infty}{{\rho_q\left(D,\delta\right)}\over {q^s}}={{
L\left(s,\chi_D\right)L\left(s,{{\delta}\over D}\right)}\over {\zeta\left
(2s\right)}}$$
for $\hbox{\rm Re$s$}>1$. We see by (1.5) that the right-hand side here
is regular for $s\neq 1$. Using also Lemma 3.4 and (1.6) we
obtain part (i) by analytic continuation.

To see part (ii) we note that $\hbox{\rm ${\rm r}{\rm e}{\rm s}_{
s=1}$}E\left(z,s\right)={3\over {\pi}}$ for every $z\in\bbb H$
by [I], (3.26). We obtain part (ii) from part (i) by analytic continuation.

To see part (iii) we apply the $D=dd'<0$ case of Theorem 1.4 of [I-L-T]. Note
that the normnalization of Fourier coefficients is
different in that paper than in the present paper,
compare [I-L-T, (1.9)] to our formulas (1.2) and (1.3). We
see in this way that our $b_j\left(n\right)$ corresponds to
$b_{\psi}\left(n\right)\left(4\pi |n|\right)^{-1/4}$ in the notation of [I-L-T]. It is also
important, as was mentioned already in the proof of
Lemma 3.4, that $\left({D\over {-1}}\right)=1$, and only the equivalence classes of
positive definite quadratic forms are considered in
[I-L-T], while we consider both positive definite and negative definite
forms. See the second paragraph above Definition 1.2 in
[I-L-T] and our formula (1.8). Finally, applying Lemma 3.2
(ii) we see that in the case $b_j\left(D\right)\neq 0$ the only $
\psi$ which
is present in the summation in [I-L-T, (1.14)] is a
constant multiple of $F_j$. In the case $b_j\left(D\right)=0$ the
summation in [I-L-T, (1.14)] is empty, and the right-hand
side of (3.14) is $0$, as needed. Taking into account these
considerations we get part (iii). The lemma is proved.

{\bf 4. Proof of Theorem 1.1. }
\medskip

{\bf 4.1. A special case.} We say that a function $\chi$ satisfies Condition $
D$ if $\chi$ is
an even entire function satisfying that for every fixed $A,B>0$
the function $\left|\chi (z)\right|e^{\left|z\right|A}$ is bounded on the strip $\left
|\hbox{\rm Im}\,z\right|\le B$.

We first prove Theorem 1.1 for such functions.

If $f$ is an automorphic function and the following
integral is absolutely convergent, define
$$\zeta (f,r):=\int_{{\cal F}_1}f(z)\overline {E\left(z,{1\over 2}
+ir\right)}d\mu_z,$$
where $E\left(z,s\right)$ is the Eisenstein series for $PSL_2({\bf Z}
)$. Let $\left\{U_l(z):\;l\ge 0\right\}$ be a complete orthonormal system of
Maass forms for $PSL_2({\bf Z})$. The function $U_0(z)$ is constant,
and $U_l(z)$ is a cusp form for $l\ge 1$. We assume that every $U_
l$ is a simultaneous Hecke
eigenform. Then by [I-K, Theorem 15.5]
we have that if $f_1$ and $f_2$ are bounded functions on ${\cal F}_
1$, then
$$\left(f_1,f_2\right)_1=\sum_{l=0}^{\infty}\left(f_1,U_l\right)_
1\overline {\left(f_2,U_l\right)_1}+{1\over {4\pi}}\int_{-\infty}^{
\infty}\zeta (f_1,\rho )\overline {\zeta (f_2,\rho )}d\rho .\eqno
(4.1)$$
We use the notations of Theorem 1.1. For $i=1,2$ let us
choose integers $n_i,t_i$ such that $n_i>0$ and $t_i^2-4n_i=\delta_
i$. Let $m_1,m_2\in {\cal K}_E$ with a large enough
absolute constant $E>0$. Assume that $m_2$ is real. Then we apply (4.1) for the functions
$$f_1(z)=M_{t_1,n_1,D_1,m_1}(z),\quad f_2(z)=M_{t_2,n_2,D_2,m_2}(
z).\eqno (4.2)$$
We then see from the last sentence of Lemma 2.6 that the contribution of
the odd cusp forms $U_l$ in (4.1) is $0$. We also see by
Lemma 3.2 that for the even cusp forms $U_l$ we can take
the functions
$${{\hbox{\rm Shim}F_j}\over {\sqrt {\left(\hbox{\rm Shim}F_j,\hbox{\rm Shim}
F_j\right)_1}}}$$
for $j\ge 1$. We see by (3.3) that
$$\Delta_0\left(\hbox{\rm Shim}F_j\right)=\left(-{1\over 4}-4r_j^
2\right)\hbox{\rm Shim}F_j.$$
On the other hand $U_0\left(z\right)=\left({3\over {\pi}}\right)^{
1/2}$  for every
$z\in\bbb H$ by [I], (3.26) and (6.33). Introduce the notations
$$A_{m_1,\delta_1,n_1}\left(\lambda\right):=\int_0^{\infty}m_1\left
({{\left|\delta_1\right|}\over {4n_1}}\sinh^2r\right)g_{\hbox{\rm }
\lambda}(r)\sinh rdr,\eqno (4.3)$$
$$A_{m_2,\delta_2,n_2}\left(\lambda\right):=\int_0^{\infty}m_2\left
({{\left|\delta_2\right|}\over {4n_2}}\sinh^2r\right)g_{\hbox{\rm }
\lambda}(r)\sinh rdr,\eqno (4.4)$$
$$H\left(\lambda\right)=H_{m_1,m_2,\delta_1,\delta_2,n_1,n_2}\left
(\lambda\right):=A_{m_1,\delta_1,n_1}\left(\lambda\right)A_{m_2,\delta_
2,n_2}\left(\lambda\right).\eqno (4.5)$$
Then we get from Lemma 2.2, (4.1), (4.2), Lemma 2.6, Lemma 3.5 and (1.7) that the sum of
$${{48}\over {\pi}}\delta_{1,D_1}\delta_{1,D_2}|\delta_1\delta_2|^{
1/2}L\left(1,\delta_1\right)L\left(1,\delta_2\right)H\left(0\right
),\eqno (4.6)$$
$$576\pi^2\sum_{j=1}^{\infty}\left(\hbox{\rm Shim}F_j,\hbox{\rm Shim}
F_j\right)_1|\delta_1\delta_2|^{3/4}b_j\left(D_1\right)\overline {
b_j\left({{\delta_1}\over {D_1}}\right)}b_j\left({{\delta_2}\over {
D_2}}\right)\overline {b_j\left(D_2\right)}H\left(-{1\over 4}-4r_
j^2\right)\eqno (4.7)$$
and
$$2\pi\int_{-\infty}^{\infty}{{L^{\ast}\left({1\over 2}-i\rho ,D_
1\right)L^{\ast}\left({1\over 2}-i\rho ,{{\delta_1}\over {D_1}}\right
)L^{\ast}\left({1\over 2}+i\rho ,D_2\right)L^{\ast}\left({1\over
2}+i\rho ,{{\delta_2}\over {D_2}}\right)H\left(-{1\over 4}-\rho^2\right
)}\over {\zeta\left(1+2i\rho\right)\zeta\left(1-2i\rho\right)}}d\rho\eqno
(4.8)$$
equals the sum of
$$4\pi E_{\delta_1,\delta_2,D_1,D_2}\int_0^{\infty}m_1\left({{\left
|\delta_1\right|}\over {n_1}}r\left(1+r\right)\right)m_2\left({{\left
|\delta_2\right|}\over {n_2}}r\left(1+r\right)\right)dr\eqno (4.9
)$$
and
$$8\sum_{f\in {\bf Z},f^2>\left|\delta_1\delta_2\right|}h_{D_1,D_
2}\left(\delta_1,\delta_2,f\right){\cal L}\left({{t_1}\over {\sqrt {
n_1}}},{{t_2}\over {\sqrt {n_2}}},\left|{f\over {\sqrt {\left|\delta_
1\delta_2\right|}}}\right|,m_1,m_2\right).\eqno (4.10)$$
By (2.33) and [G-R], p 999, 9.133 we have for $\lambda =-{1\over
4}-\tau^2$ that
$$g_{\lambda}(r)=F\left({1\over 4}+{{i\tau}\over 2},{1\over 4}-{{
i\tau}\over 2},1;-\sinh^2r\right)\eqno $$
for every $\lambda\le 0$ and $r\ge 0$. Making the substitution
$x=\sinh^2r$ we then get by (4.3) and (4.4) that
$$A_{m_1,\delta_1,n_1}\left(\lambda\right)=\int_0^{\infty}m_1\left
({{\left|\delta_1\right|}\over {4n_1}}x\right)F\left({1\over 4}+{{
i\tau}\over 2},{1\over 4}-{{i\tau}\over 2},1;-x\right){{dx}\over {
2\sqrt {1+x}}},\eqno (4.11)$$
$$A_{m_2,\delta_2,n_2}\left(\lambda\right)=\int_0^{\infty}m_2\left
({{\left|\delta_2\right|}\over {4n_2}}x\right)F\left({1\over 4}+{{
i\tau}\over 2},{1\over 4}-{{i\tau}\over 2},1;-x\right){{dx}\over {
2\sqrt {1+x}}}.$$
If we fix $C$ to be a large enough absolute constant, then
we can choose
$$\hbox{\rm \ $m_2\left(y\right)=$$\left(1+{{4n_2}\over {\left|\delta_
2\right|}}y\right)^{-C}$},\eqno (4.12)$$
since then $m_2\in {\cal K}_E$. Then we have
$$A_{m_2,\delta_2,n_2}\left(\lambda\right)={{\Gamma\left(C-{1\over
4}\pm{{i\tau}\over 2}\right)}\over {2\Gamma\left(C\right)\Gamma\left
(C+{1\over 2}\right)}}\eqno (4.13)$$
by [G-R], p 807, 7.512.10.

Let $\chi$ be a given function satisfying Condition $D$. Let us
choose $m_1$ such that
$$m_1\left({{\left|\delta_1\right|}\over {4n_1}}x\right){1\over {
2\sqrt {1+x}\Gamma\left(C\right)\Gamma\left(C+{1\over 2}\right)}}\eqno
(4.14)$$
equals
$${1\over {\pi}}\int_0^{\infty}F\left({3\over 4}-iz,{3\over 4}+iz
,1,-x\right)\left|{{\Gamma\left({1\over 4}+iz\right)\Gamma\left({
3\over 4}+iz\right)}\over {\Gamma\left(2iz\right)}}\right|^2{{\chi\left
(z\right)}\over {\Gamma\left(C-{1\over 4}\pm iz\right)}}dz\eqno (
4.15)$$
for every $x\ge 0$. The function ${{\chi\left(z\right)}\over {\Gamma\left
(C-{1\over 4}\pm iz\right)}}$ also satisfies
Condition $D$. It follows then by Lemma 3.7 of [B5], by [G-R], p 998, 9.131.1
and by (4.11) that $m_1\in {\cal K}_E$ and
$$A_{m_1,\delta_1,n_1}\left(\lambda\right)=\chi\left({{\tau}\over
2}\right){{2\Gamma\left(C\right)\Gamma\left(C+{1\over 2}\right)}\over {
\Gamma\left(C-{1\over 4}\pm{{i\tau}\over 2}\right)}}.$$
Then by (4.5) and (4.13) we get for $\lambda =-{1\over 4}-\tau^2$ that
$$H\left(\lambda\right)=\chi\left({{\tau}\over 2}\right).\eqno (4
.16)$$
We now examine the function  $\!{\cal L}\left(\tau_1,\tau_2,\phi
,m_1,m_2\right)$ defined in
(2.7) and (2.8). We note that
$$2\Phi\left(2r_1+1\right)\left(2r_2+1\right)-\Phi^2-\left(2r_1+1\right
)^2-\left(2r_2+1\right)^2+1$$
equals
$$\left(2r_2+1-a\left(r_1,\Phi\right)\right)\left(b\left(r_1,\Phi\right
)-2r_2-1\right),$$
where
$$a\left(r_1,\Phi\right):=\left(2r_1+1\right)\Phi -\sqrt {\left(\Phi^
2-1\right)\left(\left(2r_1+1\right)^2-1\right)},\eqno (4.17)$$
$$b\left(r_1,\Phi\right):=\left(2r_1+1\right)\Phi +\sqrt {\left(\Phi^
2-1\right)\left(\left(2r_1+1\right)^2-1\right)}.\eqno (4.18)$$
Then we have that the set (2.8) can be written as
$$\left\{\left(r_1,r_2\right)\in {\bf R}_{+}^2:\;{{a\left(r_1,\Phi\right
)-1}\over 2}\le r_2\le{{b\left(r_1,\Phi\right)-1}\over 2}\right\}
.$$
Then for the functions $m_1$ and $m_2$ defined in (4.14), (4.15)
and (4.12) we have that
$${\cal L}\left({{t_1}\over {\sqrt {n_1}}},{{t_2}\over {\sqrt {n_
2}}},\left|{f\over {\sqrt {\left|\delta_1\delta_2\right|}}}\right
|,m_1,m_2\right)\eqno (4.19)$$
equals
$$\int_0^{\infty}m_1\left({{\left|\delta_1\right|}\over {n_1}}r_1\left
(r_1+1\right)\right)\int_{{{a\left(r_1,\Phi\right)-1}\over 2}}^{{{
b\left(r_1,\Phi\right)-1}\over 2}}{{\left(1+2r_2\right)^{-2C}dr_2}\over \sqrt { {\left
(2r_2+1-a\left(r_1,\Phi\right)\right)\left(b\left(r_1,\Phi\right)
-2r_2-1\right)}}}dr_1\eqno (4.20)$$

with the notation
$$\phi :=\left|{f\over {\sqrt {\left|\delta_1\delta_2\right|}}}\right
|.\eqno (4.21)$$
In the inner integral in (4.20) we use the substitution
$q={{2r_2+1-a\left(r_1,\Phi\right)}\over {b\left(r_1,\Phi\right)-
a\left(r_1,\Phi\right)}}$, and we get that the inner integral
equals
$${1\over 2}\int_0^1{{\left(a\left(r_1,\Phi\right)+q\left(b\left(
r_1,\Phi\right)-a\left(r_1,\Phi\right)\right)\right)^{-2C}}\over {\sqrt {
q\left(1-q\right)}}}dq.$$
By [G-R], p 995, 9.111 this equals
$${{\Gamma^2\left({1\over 2}\right)}\over 2}a\left(r_1,\Phi\right
)^{-2C}F\left({1\over 2},2C,1;-{{b\left(r_1,\Phi\right)-a\left(r_
1,\Phi\right)}\over {a\left(r_1,\Phi\right)}}\right),$$
and then applying [G-R], p 999, 9.134.1 and (4.17), (4.18) we
finally get that the inner integral in (4.20) equals
$${{\Gamma^2\left({1\over 2}\right)}\over 2}\left(\left(2r_1+1\right
)\Phi\right)^{-2C}F\left(C+{1\over 2},C,1;{{\left(\Phi^2-1\right)\left
(\left(2r_1+1\right)^2-1\right)}\over {\left(2r_1+1\right)^2\Phi^
2}}\right).$$
Then applying (4.20), (4.14) and (4.15) we obtain that (4.19)
equals
$$\Phi^{-2C}\Gamma\left(C\right)\Gamma\left(C+{1\over 2}\right)\int_
0^{\infty}\left|{{\Gamma\left({1\over 4}+iz\right)\Gamma\left({3\over
4}+iz\right)}\over {\Gamma\left(2iz\right)}}\right|^2{{\chi\left(
z\right)}\over {\Gamma\left(C-{1\over 4}\pm iz\right)}}I\left(z\right
)dz\eqno (4.22)$$
with the abbreviation
$$I\left(z\right):=\int_0^{\infty}{{F\left({3\over 4}-iz,{3\over
4}+iz,1,-4r_1\left(r_1+1\right)\right)F\left(C+{1\over 2},C,1;{{\left
(\Phi^2-1\right)\left(\left(2r_1+1\right)^2-1\right)}\over {\left
(2r_1+1\right)^2\Phi^2}}\right)}\over {\left(2r_1+1\right)^{2C-1}}}
dr_1.$$
We make the substitution $x=4r_1\left(r_1+1\right)$. Then using also
[G-R], p 998, 9.131.1 we get that
$$I\left(z\right)={1\over 4}\int_0^{\infty}{{F\left({1\over 4}-iz
,{1\over 4}+iz,1,-x\right)F\left(C+{1\over 2},C,1;{{\left(\Phi^2-
1\right)x}\over {\left(x+1\right)\Phi^2}}\right)}\over {\left(x+1\right
)^{C+{1\over 2}}}}dx.\eqno (4.23)$$
We compute this integral in the following lemma. During
its proof we need the notation
$$\hbox{\rm $_3$$F_2\left(\matrix{a_1,a_2,a_3\cr
b_1,b_2\cr}
;1\right):=$}\sum_{k=0}^{\infty}{{\left(a_1\right)_k\left(a_2\right
)_k\left(a_3\right)_k}\over {n!\left(b_1\right)_k\left(b_2\right)_
k}}.$$
Here $\left(a\right)_k:={{\Gamma\left(a+k\right)}\over {\Gamma\left
(a\right)}}$ and the $b_i$ are not nonpsitive
integers. We will need only the case when one of the $a_i$ is a nonpositive
integer. In this case we have in fact a finite sum.

{\bf LEMMA 4.1.} {\it Let} $z${\it ,} $C$ {\it and} $\Phi$ {\it be real numbers such that} $
C>{1\over 4}$ {\it and} $\Phi >1${\it . Then}
$$\int_0^{\infty}F\left({1\over 4}-iz,{1\over 4}+iz,1,-x\right)F\left
(C+{1\over 2},C,1;{{\left(\Phi^2-1\right)x}\over {\Phi^2\left(1+x\right
)}}\right)\left(1+x\right)^{-C-{1\over 2}}dx\eqno (4.24)$$
{\it equals}
$${{\Gamma\left(C-{1\over 4}\pm iz\right)}\over {\Gamma\left(C\right
)\Gamma\left(C+{1\over 2}\right)}}\Phi^{2C}F\left({1\over 4}-iz,{
1\over 4}+iz,1,1-\Phi^2\right).\eqno (4.25)$$
{\it Proof.\/} We can clearly write (4.24) as the sum
$$\sum_{n=0}^{\infty}{{\left(C+{1\over 2}\right)_n\left(C\right)_
n}\over {n!n!}}a_n\left(1-{1\over {\Phi^2}}\right)^n$$
with
$$a_n:=\int_0^{\infty}F\left({1\over 4}-iz,{1\over 4}+iz,1,-x\right
)x^n\left(1+x\right)^{-C-{1\over 2}-n}dx.$$
On the other hand, we have
$$\Phi^{2C}F\left({1\over 4}-iz,{1\over 4}+iz,1,1-\Phi^2\right)=\Phi^{
2C-{1\over 2}-2iz}F\left({3\over 4}+iz,{1\over 4}+iz,1,1-{1\over {
\Phi^2}}\right)$$
by [G-R], p 998, 9.131.1, and
$$\Phi^{2C-{1\over 2}-2iz}=\sum_{r=0}^{\infty}{{\left(C-{1\over 4}
-iz\right)_r}\over {r!}}\left(1-{1\over {\Phi^2}}\right)^r,$$
since the rght-hand side here is a binomial series.
Therefore (4.25) equals
$${{\Gamma\left(C-{1\over 4}\pm iz\right)}\over {\Gamma\left(C\right
)\Gamma\left(C+{1\over 2}\right)}}\sum_{n=0}^{\infty}{{\left(C-{1\over
4}-iz\right)_n}\over {n!}}\left(1-{1\over {\Phi^2}}\right)^n\,_3F_
2\left(\matrix{-n,{1\over 4}+iz,{3\over 4}+iz\cr
1,{5\over 4}-C+iz-n\cr}
;1\right).$$
So it is enough to show that
$$a_n={{n!\Gamma\left(C-{1\over 4}+iz\right)\Gamma\left(C-{1\over
4}-iz+n\right)}\over {\Gamma\left(C+n\right)\Gamma\left(C+{1\over
2}+n\right)}}\,_3F_2\left(\matrix{-n,{1\over 4}+iz,{3\over 4}+iz\cr
1,{5\over 4}-C+iz-n\cr}
;1\right)\eqno (4.26)$$
for every $n\ge 0$. Writing
$$x^n\left(1+x\right)^{-n}=\sum_{k=0}^n{{\left(-n\right)_k}\over {
k!}}\left({1\over {1+x}}\right)^k$$
by the binomial theorem and applying [G-R], p 807, 7.512.10 we get that
$$a_n={{\Gamma\left(C-{1\over 4}\pm iz\right)}\over {\Gamma\left(
C+{1\over 2}\right)\Gamma\left(C\right)}}\,_3F_2\left(\matrix{-n,
C-{1\over 4}+iz,C-{1\over 4}-iz\cr
C+{1\over 2},C\cr}
;1\right).\eqno (4.27)$$
We have to show that the right-hand sides of (4.26) and
(4.27) are the same. Now, the right-hand side of (4.27)
equals
$${{\Gamma\left(C-{1\over 4}\pm iz\right)\Gamma\left(1+n\right)}\over {
\Gamma\left(C+{1\over 2}+n\right)\Gamma\left(C\right)}}\,_3F_2\left
(\matrix{-n,{1\over 4}+iz,{1\over 4}-iz\cr
1,C\cr}
;1\right)\eqno (4.28)$$
by Corollary 3.3.5 of [A-A-R]. We see that (4.28) equals the
right-hand sides of (4.26) by [S], p 121, (4.3.4.2). The
lemma is proved.

By (4.22), (4.23), Lemma 4.1 and (4.21) we get that (4.19) equals
$${1\over 4}\int_0^{\infty}\left|{{\Gamma\left({1\over 4}+iz\right
)\Gamma\left({3\over 4}+iz\right)}\over {\Gamma\left(2iz\right)}}\right
|^2F\left({1\over 4}-iz,{1\over 4}+iz,1,1-{{f^2}\over {\left|\delta_
1\delta_2\right|}}\right)\chi\left(z\right)dz.\eqno (4.29)$$
We now compute
$$\int_0^{\infty}m_1\left({{\left|\delta_1\right|}\over {n_1}}r\left
(1+r\right)\right)m_2\left({{\left|\delta_2\right|}\over {n_2}}r\left
(1+r\right)\right)dr.\eqno (4.30)$$
By (4.12), (4.14) and (4.15) we have that (4.30) equals
$${{2\Gamma\left(C\right)\Gamma\left(C+{1\over 2}\right)}\over {\pi}}
\int_0^{\infty}\left|{{\Gamma\left({1\over 4}+iz\right)\Gamma\left
({3\over 4}+iz\right)}\over {\Gamma\left(2iz\right)}}\right|^2{{\chi\left
(z\right)}\over {\Gamma\left(C-{1\over 4}\pm iz\right)}}J\left(z\right
)dz$$
with the abbreviation
$$J\left(z\right):=\int_0^{\infty}F\left({3\over 4}-iz,{3\over 4}
+iz,1,-4r\left(1+r\right)\right)\left(1+2r\right)^{1-2C}dr.$$
Applying the substitution $x=4r\left(r+1\right)$ we get
$$J\left(z\right)={1\over 4}\int_0^{\infty}F\left({3\over 4}-iz,{
3\over 4}+iz,1,-x\right)\left(1+x\right)^{-C}dx,$$
so applying [G-R], p 807, 7.512.10 we get that (4.30) equals
$${1\over {2\pi}}\int_0^{\infty}\left|{{\Gamma\left({1\over 4}+iz\right
)\Gamma\left({3\over 4}+iz\right)}\over {\Gamma\left(2iz\right)}}\right
|^2\chi\left(z\right)dz.\eqno (4.31)$$
Then by (4.6)-(4.10), (4.16), (4.19), (4.29), (4.30), (4.31) we
get Theorem 1.1 for $\chi$ satisfying Condition $D$.

{\bf 4.2. The end of the proof.} We now extend the theorem for the general case. For
this sake we first need the following upper bound.

{\bf LEMMA 4.2.} {\it There is an absolute constant} $C>0$ {\it such that the sequence}
$$\left|\left(\hbox{\rm Shim}F_j,\hbox{\rm Shim}F_j\right)_1b_j\left
(D_1\right)b_j\left({{\delta_1}\over {D_1}}\right)b_j\left({{\delta_
2}\over {D_2}}\right)b_j\left(D_2\right)\right|\left(1+r_j\right)^{
-C}\eqno (4.32)$$
{\it is bounded for} $j\ge 1${\it , and the sequence}
$$\left|h_{D_1,D_2}\left(\delta_1,\delta_2,f\right)\right|\left(1
+f^2\right)^{-C}\eqno (4.33)$$
{\it is bounded for} $f\in {\bf Z},f^2>\left|\delta_1\delta_2\right
|$.

{\it Proof.\/} It follows from Theorem 5 of [Du] that there is
an absolute constant $C_1>0$ such that the sequence
$$\left|b_j\left(D_1\right)b_j\left({{\delta_1}\over {D_1}}\right
)b_j\left({{\delta_2}\over {D_2}}\right)b_j\left(D_2\right)\right
|e^{-2\pi r_j}\left(1+r_j\right)^{-C_1}\eqno (4.34)$$
is bounded for $j\ge 1$.

Let us write
$$u_j:={{\hbox{\rm Shim$F_j$}}\over {\sqrt {\left(\hbox{\rm Shim$
F_j$},\hbox{\rm Shim$F_j$}\right)_1}}},$$
then $u_j$ is a Maass cusp form of weight $0$, we have
$\left(\hbox{\rm $u_j$},\hbox{\rm $u_j$}\right)_1=1$, and by (3.3) we see that $
\Delta_0u_j=\left(-{1\over 4}-4r_j^2\right)u_j$.
We clearly have
$$\left(\hbox{\rm Shim$F_j$},\hbox{\rm Shim$F_j$}\right)_1={1\over {\left
|\rho_{u_j}(1)\right|^2}},$$
see (1.2). By [I], (8.1), (8.5) and (8.43) we then get that there is
an absolute constant $C_2>0$ such that the sequence
$$\left(\hbox{\rm Shim$F_j$},\hbox{\rm Shim$F_j$}\right)_1e^{2\pi
r_j}\left(1+r_j\right)^{-C_2}\eqno (4.35)$$
is bounded for $j\ge 1$. By (4.34) and (4.35) we obtain (4.32). The estimate (4.33) follows at once from Lemma 3.1 of
[B2].

The proof of the following lemma is very similar to the proof of lemma 3.7 of [B5].

{\bf LEMMA 4.3.} {\it Let} $A>0$ {\it be given. Then there is a} $
\beta >0$ {\it depending only on} $A$ {\it such that the following statement holds. If} $
M$ {\it is a nonnegative function on} $[0,\infty )$ {\it satisfying that the function} $
M(R)\left(1+R\right)^{\beta}$ {\it is bounded on} $[0,\infty )${\it , and} $
\chi$ {\it is any even holomorphic function on the strip} $\left|\hbox{\rm Im}\,
z\right|<\beta$ $wi${\it th} $\left|\chi\left(z\right)\right|\le
M\left(\left|z\right|\right)$ {\it on this strip, then we have that}
$$T_{\chi}\left(u\right)\ll_{\beta ,M}\left(1+u\right)^{-A}$$
{\it for} $u\ge 0$.

{\it Proof.\/} By [S], (1.8.1.11) we know for real $z$ that
$$F\left({3\over 4}-iz,{3\over 4}+iz,1,-u\right)\left|{{\Gamma\left
({1\over 4}+iz\right)\Gamma\left({3\over 4}+iz\right)}\over {\Gamma\left
(2iz\right)}}\right|^2=\phi (u,z)+\phi (u,-z),$$
where
$$\phi (u,z)={{\Gamma\left({1\over 4}-iz\right)\Gamma\left({3\over
4}-iz\right)}\over {\Gamma\left(-2iz\right)}}u^{iz-{3\over 4}}F\left
({3\over 4}-iz,{3\over 4}-iz,1-2iz,-{1\over u}\right).$$
Hence using also [G-R], p 998, 9.131.1 we have that
$$T_{\chi}\left(u\right)=\left(1+u\right)^{{1\over 2}}{1\over {2\pi}}
\int_{-\infty}^{\infty}\phi (u,z)\chi\left(z\right)dz.$$
We push the line of integration upwards to a line $\hbox{\rm Im$\,
z=B$ }$with a large positive number
$B$ depending on $A$. Using  [G-R], p. 995, 9.111 to estimate
$\phi (u,z)$ we obtain the lemma.

We also need the following lemma, proved in [B5].

{\bf LEMMA 4.4.} {\it Let} $\beta >0$ {\it and let} $\chi$ {\it be an even holomorphic function on the strip} $\left
|\hbox{\rm Im}\,z\right|<\beta$ {\it such that for a fixed} $A>0$ {\it the function} $\left
|\chi (z)\right|e^{A\left|z\right|^2}$ {\it is bounded}
{\it on the strip} $\left|\hbox{\rm Im}\,z\right|<\beta$. {\it Then for every} $
0<\gamma <\beta$ {\it there is a sequence} $\chi_n$ {\it of entire functions, and a nonnegative function} $
M$ {\it on  }
$[0,\infty )$ {\it with the following properties. The function} $
\chi_n$
{\it satisfies Condition D for every} $n$, {\it for every fixed} $
K>0$
{\it the function} $M(R)e^{KR}$ {\it is bounded on} $[0,\infty )${\it , we have }
$\left|\chi_n\left(z\right)\right|\le M\left(\left|z\right|\right
)$ {\it for every} $n\ge 1$ {\it and} $\left|\hbox{\rm Im}\,z\right
|<\gamma${\it , and finally,} $\chi_n(z)\rightarrow\chi (z)$ {\it for every} $\left
|\hbox{\rm Im}\,z\right|<\gamma${\it .}

{\it Proof.\/} See [B5], Lemma 5.1.

We now finish the proof of Theorem 1.1. Our argument is
similar to that applied in [B5].

Let $\beta$ be a large enough absolute constant, and let $\chi$ be a function
satisfying Condition $A_{\beta}$. Then we easily see using Lemmas 4.2, 4.3 and the dominated convergence theorem that it is enough to prove Theorem
1.1 for every function $\chi (z)e^{-z^2/N}$ ($N$ is a positive integer) instead of $
\chi$. So we may assume that there is an $A>0$ such that $\chi (z
)e^{A\left|z\right|^2}$
is bounded on the strip $\left|\hbox{\rm Im$\,z$}\right|<\beta$. Finally, for such
functions the theorem follows from Lemmas 4.4, 4.3, 4.2, the dominated convergence theorem and the already
proved special case of Theorem 1.1. The theorem is proved.

{\bf 5. Proof of Theorem 1.2. }
\medskip

Recall the notations from Subsection 1.8. It is easy to see that for any $T\in SL_2({\bf R})$ we have
$${{h\left(Tz,Tw\right)}\over {h\left(z,w\right)}}=\left({{j_T(w)}\over {\left
|j_T(w)\right|}}\right)^2\left({{j_T(z)}\over {\left|j_T(z)\right
|}}\right)^{-2}.\eqno (5.1)$$
One can see easily using (5.1) and [I], (1.10) that if $\tau\in SL_
2({\bf R})$, then
$$h\left(\tau^{-1}\gamma\tau z,z\right)\left({{j_{\tau^{-1}\gamma
\tau}(z)}\over {\left|j_{\tau^{-1}\gamma\tau}(z)\right|}}\right)^
2=h\left(\gamma\tau z,\tau z\right)\left({{j_{\gamma}(\tau z)}\over {\left
|j_{\gamma}(\tau z)\right|}}\right)^2.\eqno (5.2)$$
Hence $N_{t,n,D,m}(z)$ is $SL_2({\bf Z})$-invariant.

As in the proof of Lemma 2.6, let

$$[\gamma ]=\left\{\tau^{-1}\gamma\tau :\,\tau\in SL_2({\bf Z})\right
\}$$
and for $\gamma\in\Gamma_{n,t}$, write
$$T_{\gamma}=\sum_{\delta\in [\gamma ]}\int_{{\cal F}_1}m\left(z,
\delta z\right)h\left(\delta z,z\right)\left({{j_{\delta}(z)}\over {\left
|j_{\delta}(z)\right|}}\right)^2u(z)d\mu_z.$$
Then we have
$$T_{\gamma}=\int_{C(\gamma )\setminus\bbb H}m\left(z,\gamma z\right
)h\left(\gamma z,z\right)\left({{j_{\gamma}(z)}\over {\left|j_{\gamma}
(z)\right|}}\right)^2u(z)d\mu_z,\eqno (5.3)$$
where
$$C\left(\gamma\right):=\left\{\tau\in SL_2({\bf Z}):\;\;\tau\gamma
=\gamma\tau\right\}.$$
It is proved on pp 117-118 of [B1] that the image of $C\left(\gamma\right
)$
in $PSL_2({\bf Z})$ is trivial if $\delta =t^2-4n$ is a square, and it is infinite
cyclic if $\delta$ is not a square.

As in the proof of [B1, Lemma 2] let $h=h_{\gamma}\in SL_2({\bf R}
)$ be such
that $h^{-1}\gamma hz=Rz$ for every $z\in\bbb H$ with an $R>1$. We then have
$$\sqrt R+{1\over {\sqrt R}}={t\over {\sqrt n}},\,\sqrt R-{1\over {\sqrt
R}}={{\sqrt {\delta}}\over {\sqrt n}},\,R+{1\over R}-2={{\delta}\over
n}.\eqno (5.4)$$

We will need later the concrete form of $h$. Let
$\gamma =\left(\matrix{a&b\cr
c&d\cr}
\right)$. If $c\neq 0$, then the two fixed points of $\gamma$ are
$$z_1:={{a-d+\sqrt {\delta}}\over {2c}},\;z_2:={{a-d-\sqrt {\delta}}\over {
2c}}.\eqno (5.5)$$
Then one can take $h=\left(\matrix{z_1&{{z_2}\over {z_1-z_2}}\cr
1&{1\over {z_1-z_2}}\cr}
\right)$, and we have
$$h\left(\infty\right)=z_1,\;h\left(0\right)=z_2.\;\eqno (5.6)$$
If $c=0$, then the two fixed points are
$$z_1:=\infty ,\;z_2:={b\over {d-a}}.$$
Then one can take $h=\left(\matrix{1&{b\over {d-a}}\cr
0&1\cr}
\right)$ if $a>d$, and $h=\left(\matrix{{b\over {d-a}}&-1\cr
1&0\cr}
\right)$
if $d>a$. So in this case we have
$$h\left(\infty\right)=z_1,\;h\left(0\right)=z_2\;\eqno (5.7)$$
if $a>d$, and
$$h\left(\infty\right)=z_2,\;h\left(0\right)=z_1\;\eqno (5.8)$$
if $d>a$.

Then by (5.2) and (5.3) we get that
$$T_{\gamma}=\int_{h^{-1}C(\gamma )h\setminus\bbb H}m\left(z,Rz\right
)h\left(Rz,z\right)u\left(hz\right)d\mu_z.$$
In case $\delta$ is not a square, let $r_0>1$ be such that
$\left(\matrix{\sqrt {r_0}&0\cr
0&1/\sqrt {r_0}\cr}
\right)$ is a generator of the image of $h^{-1}C(\gamma )h$ in
$PSL_2({\bf R})$. Let $I_{\gamma}=[1,r_0)$ if $\delta$ is not a square, and let
$I_{\gamma}=(0,\infty )$ otherwise. Then by the substitution
$$z=re^{i\left({{\pi}\over 2}+\theta\right)}\eqno (5.9)$$
we have that
$$T_{\gamma}=\int_{-\pi /2}^{\pi /2}\int_{I_{\gamma}}m\left({{\delta}\over {
4n\cos^2\theta}}\right)h\left(Rz,z\right)u\left(h\left(re^{i\left
({{\pi}\over 2}+\theta\right)}\right)\right){{drd\theta}\over {r\cos^
2\theta}},$$
where $z$ is given by (5.9). Now, by (1.19) and (5.4) we see that
$$h\left(Rre^{i\left({{\pi}\over 2}+\theta\right)},re^{i\left({{\pi}\over
2}+\theta\right)}\right)={{\left(-{{\sqrt {\delta}}\over t}\sin\theta
+i\cos\theta\right)^2}\over {\left|-{{\sqrt {\delta}}\over t}\sin
\theta +i\cos\theta\right|^2}}.$$
Hence we have
$$T_{\gamma}=\int_{-\pi /2}^{\pi /2}m\left({{\delta}\over {4n\cos^
2\theta}}\right){{\left(-{{\sqrt {\delta}}\over t}\sin\theta +i\cos
\theta\right)^2}\over {\left|-{{\sqrt {\delta}}\over t}\sin\theta
+i\cos\theta\right|^2}}F_{\gamma}\left(e^{i\left({{\pi}\over 2}+\theta\right
)}\right){{d\theta}\over {\cos^2\theta}}$$
with
$$F_{\gamma}\left(z\right):=\int_{I_{\gamma}}u\left(h\left(rz\right
)\right){{dr}\over r}\eqno (5.10)$$
for $z\in\bbb H$. Hence we proved that
$$\int_{{\cal F}_1}N_{t,n,D,m}(z)u(z)d\mu_z\eqno (5.11)$$
equals
$$\sum_{[\gamma ]}\omega_D\left(\gamma\right)\int_{-\pi /2}^{\pi
/2}m\left({{\delta}\over {4n\cos^2\theta}}\right){{\left(-{{\sqrt {
\delta}}\over t}\sin\theta +i\cos\theta\right)^2}\over {\left|-{{\sqrt {
\delta}}\over t}\sin\theta +i\cos\theta\right|^2}}F_{\gamma}\left
(e^{i\left({{\pi}\over 2}+\theta\right)}\right){{d\theta}\over {\cos^
2\theta}},\eqno (5.12)$$
where the summation is over the $SL_2({\bf Z})$-conjugacy classesof $
\Gamma_{n,t}$.

Let $f_{\lambda}(\theta )$ be the unique even solution of the equation (1.20)
with $f_{\lambda}(0)=1$. It is proved on p 119 of [B1] with slightly different
notations that
$$F_{\gamma}\left(e^{i\left({{\pi}\over 2}+\theta\right)}\right)=
F_{\gamma}\left(e^{i{{\pi}\over 2}}\right)f_{\lambda}(\theta )+\left
({d\over {d\theta}}\left(F_{\gamma}\left(e^{i\left({{\pi}\over 2}
+\theta\right)}\right)\right)\right)\left(0\right)h_{\lambda}(\theta
).\eqno (5.13)$$
It is clear that
$$F_{\gamma}\left(e^{i{{\pi}\over 2}}\right)=\int_{C_{Q_{\gamma}}}
udS,\eqno (5.14)$$
where $dS={{\left|dz\right|}\over y}$ is the hyperbolic arc length.

For $z\in\bbb H$ define
$$u_h\left(z\right):=u\left(hz\right).$$
Then by (5.10) we have
$$i\left({d\over {d\theta}}\left(F_{\gamma}\left(e^{i\left({{\pi}\over
2}+\theta\right)}\right)\right)\right)\left(0\right)=-\int_{I_{\gamma}}{{
\partial u_h}\over {\partial x}}\left(z\right)dz.$$
Since $u_h\left(z\right)$ takes the same values at the endpoints of
$I_{\gamma}$, we can write also
$$i\left({d\over {d\theta}}\left(F_{\gamma}\left(e^{i\left({{\pi}\over
2}+\theta\right)}\right)\right)\right)\left(0\right)=-2\int_{I_{\gamma}}{{
\partial u_h}\over {\partial z}}\left(z\right)dz,$$
where we write
$${{\partial}\over {\partial z}}={1\over 2}\left({{\partial}\over {
\partial x}}-i{{\partial}\over {\partial y}}\right).$$
Using (5.5)-(5.8) we then get that
$$i\left({d\over {d\theta}}\left(F_{\gamma}\left(e^{i\left({{\pi}\over
2}+\theta\right)}\right)\right)\right)\left(0\right)=2\int_{C_{Q_{
\gamma}}}{{\partial u}\over {\partial z}}\left(z\right)dz.\eqno (5.15)$$
By (5.11)-(5.15) and by the remarks in Subsection 1.6 on the correspondence between $
\Gamma_{n,t}$
and $\scr{Q}_{\delta}$ we obtain that (5.11) equals
$$\left(\sum_{Q\in\Lambda_{\delta}}\omega_D\left(Q\right)\int_{C_
Q}udS\right)F_1\left(\lambda\right)+\left(\sum_{Q\in\Lambda_{\delta}}
\omega_D\left(Q\right)\int_{C_Q}{{\partial u}\over {\partial z}}\left
(z\right)dz\right)F_2\left(\lambda\right)\eqno (5.16)$$
with
$$F_1\left(\lambda\right):=\int_{-\pi /2}^{\pi /2}m\left({{\delta}\over {
4n\cos^2\theta}}\right){{\left(-{{\sqrt {\delta}}\over t}\sin\theta
+i\cos\theta\right)^2}\over {\left|-{{\sqrt {\delta}}\over t}\sin
\theta +i\cos\theta\right|^2}}f_{\lambda}(\theta ){{d\theta}\over {\cos^
2\theta}},$$
$$F_2\left(\lambda\right):=-2i\int_{-\pi /2}^{\pi /2}m\left({{\delta}\over {
4n\cos^2\theta}}\right){{\left(-{{\sqrt {\delta}}\over t}\sin\theta
+i\cos\theta\right)^2}\over {\left|-{{\sqrt {\delta}}\over t}\sin
\theta +i\cos\theta\right|^2}}h_{\lambda}(\theta ){{d\theta}\over {\cos^
2\theta}}.$$
We now show that
$$\sum_{Q\in\Lambda_{\delta}}\omega_D\left(Q\right)\int_{C_Q}udS=
0.\eqno (5.17)$$
Indeed, since $D<0$, we have $\left({D\over {-1}}\right)=-1$, see [D], p 41.
Therefore $\omega_D\left(Q\right)=-\omega_D\left(-Q\right)$. But $
\int_{C_Q}udS=\int_{C_{-Q}}udS$,
because we integrate here with respect to the arc
length, so the orientation of the curves is not relevant.
Hence (5.17) follows.

Taking into account that $h_{\lambda} (\theta )$ is odd and $
t=\sqrt {\delta +4n}$ one can compute that
$$F_2\left(\lambda\right)=-4\int_{-\pi /2}^{\pi /2}m\left({{\delta}\over {
4n\cos^2\theta}}\right){{\sqrt {1+{{4n}\over {\delta}}}}\over {1+{{
4n}\over {\delta}}\cos^2\theta}}h_{\lambda}(\theta ){{\sin\theta
d\theta}\over {\cos\theta}}.\eqno (5.18)$$
We have
$${1\over {\left(u,u\right)_1}}\sum_{Q\in\Lambda_{\delta}}\omega_
D\left(Q\right)\int_{C_Q}{{\partial u}\over {\partial z}}\left(z\right
)dz={1\over i}12\sqrt {\pi}\delta ^{3/4}\overline {b_j\left(D\right
)}b_j\left({{\delta}\over D}\right).\eqno (5.19)$$
Indeed, this is proved in Proposition 6 of [D-I-T] and Theorem 1.4 of [I-L-T].

By (5.11), (5.16), (5.17), (5.18) and (5.19) we get the theorem.

 \bigskip\noindent {\bf References}

\nobreak
\parindent=12pt
\nobreak

\item{[A-A]} S. Ahlgren, N. Andersen, {\it Kloosterman sums and Maass cusp forms of half integral weight for the modular group}, International Mathematics Research Notices 2018.2 (2018), 492-570.

\item{[A-A-R]} G.E. Andrews, R. Askey, R. Roy, {Special
Functions}, {\it Cambridge Univ. Press,\/} 1999

\item{[A-D]} N. Andersen, W. Duke, {\it Modular invariants for real quadratic fields and Kloosterman sums}, Alg. Number Theory 14.6 (2020), 1537-1575.

\item{[B1]} A. Bir\'o, {\it Cycle integrals of Maass forms of weight 0 and Fourier coefficients of Maass forms of weight 1/2}, Acta Arithmetica,
94 (2) (2000), 103-152.

\item{[B2]} A. Bir\'o, {\it Local square mean in the hyperbolic circle problem}, arXiv e-prints (2024), arXiv:2403.16113v2, to appear in Algebra and Number Theory

\item{[B3]} A. Bir\'o, {\it On a generalization of the Selberg trace formula}, Acta Arithmetica, 87 (2) (1999), 319-338.

\item{[B4]} A. Bir\'o,  {\it A relation between triple products of weight 0 and weight 1/2 cusp forms}, Israel J. of Math., 182
(2011), 61-101.

\item{[B5]} A. Bir\'o,  {\it Triple product integrals and Rankin-Selberg
$L$-functions}, Funct. Approx. Comment. Math.  72 (2) (2025), 225-275.

\item{[B6]} A. Bir\'o, {\it On the class number of pairs of binary quadratic forms}, preprint, 2024, to appear in J. de Th. des Nombres de Bordeaux

\item{[B-C]} V. Blomer, A. Corbett, {\it A symplectic restriction problem}, Math. Annalen 382 (2022), 1323-1424.

\item{[B-M]} E.M. Baruch, Z. Mao, {\it A generalized Kohnen-Zagier formula for Maass forms}, J. London Math. Soc. (2), 82 (2010), no. 1, 1-16.

\item{[D]} H. Davenport, {Multiplicative Number Theory, Second edition,} {\it Springer}, 1980

\item{[D-I-T]} W. Duke, O. Imamoglu, \'A. T\'oth,
{\it Geometric invariants for real
quadratic fields}, Ann. of Math., 184 no. 3, (2016),
949-990.

\item{[Du]} W. Duke, {\it Hyperbolic distribution problems
and half-integral weight Maass forms}, Invent Math.,
92, 73-90 (1988)

\item{[G-R]} I.S. Gradshteyn, I.M. Ryzhik, {Table of integrals, series and
products, 6th edition,} {\it Academic Press}, 2000

\item{[H]} D.A. Hejhal, {The Selberg Trace Formula for $PSL(2,{\bf R}
)$,
vol. 2}, Springer, 1983

\item{[I]} H. Iwaniec, {Introduction to the spectral theory of automorphic forms,} {\it Rev. Mat. Iberoamericana}, 1995

\item{[I-K]} H. Iwaniec, E. Kowalski, {Analytic Number Theory,} {\it AMS Colloqium Publications, Vol. 53, Providence RI, American Mathematical Society},
2004

\item{[I-L-T]} O. Imamoglu, A Lageler, \'A. T\'oth {\it  }
{\it The Katok-Sarnak formula for higher weights}, J. of Number Theory., 235, (2022),
242-274.

\item{[K1]} W. Kohnen, {\it Fourier coefficients of modular forms of half-integral weight,} Math. Ann., 271 (1985), 237-268.

\item{[K2]} W. Kohnen, {\it Newforms of half-integral weight,\/} J. Reine Angew. Math., 333 (1982), 32-72.

\item{[Ko]} T.H. Koornwinder, {\it A
new proof of a Paley-Wiener type theorem for the Jacobi
transform,} Ark. Mat., 13 (1975), 145-159.

\item{[K-S]} S. Katok, P. Sarnak, {\it Heegner points, cycles and Maass forms,\/} Israel J. of Math., 84 (1993), 193-227.

\item{[M]} J. Morales, {\it The classification of pairs of binary quadratic forms}, Acta Arithmetica. 59 (2) (1991), 105-121.

\item{[P]} N.V. Proskurin, {\it On general Klosterman
sums (in Russian)}, Zap. Naucn. Sem. LOMI 302 (2003),
107-134.

\item{[S]} L.J. Slater, {Generalized hypergeometric functions,} {\it Cambridge Univ. Press},
1966

\item{[S-Y]} K. Soundararajan, M.P. Young {\it The prime geodesic theorem,} J. Reine Angew. Math., 676 (2013), 105-120.

\bye